\documentclass[11pt,a4paper]{amsart}
\usepackage[margin=2.5cm]{geometry}
\usepackage{amssymb,mathtools}
\usepackage{enumitem}
\usepackage{hyperref}
\usepackage{microtype}

\newtheorem{theorem}{Theorem}

\newtheorem{lemma}{Lemma}

\newtheorem{remark}{Remark}

\title[The sharp $L^p$ bound for the imaginary part of the Beurling--Ahlfors operator]
{The sharp $L^p$ bound for the imaginary part of the Beurling--Ahlfors operator on the lattice}
\author[K.~Domelevo]{Komla Domelevo}
\address{Department of Mathematics, University of W\"{u}rzburg, W\"{u}rzburg, Germany}
\email{komla.domelevo@uni-wuerzburg.de \textrm{(K. Domelevo)}}
\author[S.~Petermichl]{Stefanie Petermichl}
\address{Department of Mathematics, University of W\"{u}rzburg, W\"{u}rzburg, Germany}
\email{stefanie.petermichl@uni-wuerzburg.de \textrm{(S. Petermichl)}}
\thanks{S.P. Partially supported by the Alexander von Humboldt foundation}

\begin{document}

\begin{abstract}
	We prove the sharp~$\ell^p$ bound~$p^\star-1$ for certain second-order discrete Riesz transforms on products of integers, in all dimensions. The family includes the imaginary part of the discrete Beurling--Ahlfors operator in dimension two as a special case. While these operators admit a martingale representation via compound Poisson jump processes and heat extensions, the martingale attached to a function and the martingale attached to its transform do not enjoy strong differential subordination, so the operator-norm estimate is not accessible through prior work.
\end{abstract}

\maketitle

\section{Introduction}
\label{sec:intro}

Sharp~$L^p$ inequalities for pairs of differentially subordinate martingales
date back to the work of Burkholder~\cite{Burkholder1984}.
There the optimal constant~$p^\star-1$ is obtained from his celebrated Burkholder functional; see also~\cite{Burkholder1979,Burkholder1991}. The estimate is very general: it holds under minimal assumptions and for continuous-time filtrations. In the presence of jumps, however, the subordination condition is restrictive. Each jump must be controlled separately, so the square bracket rather than the more forgiving predictable quadratic variation is the relevant quantity.

The relation between differentially subordinate martingales and Calderón--Zygmund operators has been known since Gundy--Varopoulos~\cite{GV1976}. Bañuelos--Wang~\cite{BW1995} were the first to exploit this connection to prove new sharp or dimension-free inequalities for singular integrals. Subsequent sharp~$L^p$ arguments rely on special functions of Pichorides~\cite{Pichorides1972} and Essén~\cite{Essen1984} when orthogonality is present in addition to differential subordination, and on Burkholder~\cite{Burkholder1979,Burkholder1984,Burkholder1991} and Wang~\cite{Wang1995} when differential subordination is the only hypothesis. Deterministic proofs of sharp~$L^p$ estimates for Calderón--Zygmund operators that invoke Burkholder's theorems are available; the Bellman-function technique was used by Nazarov--Volberg~\cite{NV2004} for certain second-order Riesz transforms in the Euclidean plane, and by Domelevo--Petermichl~\cite{DP2014} in the corresponding setting on discrete abelian groups. Given how restrictive subordination becomes in the presence of jumps, it is not surprising that further difficulties arise when the operators act on functions defined on the integers.

The continuous world is over-determined: a single operator (Hilbert transform, Riesz transform, Beurling--Ahlfors transform) simultaneously satisfies many strong properties---a specific Fourier multiplier, a singular integral kernel, a martingale representation, algebraic intertwining with~$\partial$ and~$\bar\partial$, Cauchy--Riemann relations, an~$L^2$-isometry, and so on. On a lattice those properties become independent and often incompatible, and one is forced to choose which subset to preserve. There are, for instance, numerous discrete Hilbert transforms: the naïve truncated-kernel definition, the Fourier-multiplier signum, and various versions encoding Cauchy--Riemann relations. Discrete Beurling--Ahlfors operators split in the same way. The singular-sum version, the probabilistic Riesz-combination version, and the exterior-calculus signature version are three different operators, and none of them inherits \emph{all} continuous properties at once.

As is typical for singular integrals on the lattice, different continuous properties therefore lead to distinct discrete realisations. The operator studied here is the one canonically associated with the probabilistic discrete second-order Riesz transforms: it is the discrete Beurling--Ahlfors transform in the martingale sense, built from finite differences. The graph Laplacian together with right and left differences produce left and right Riesz transforms~$R_k^\pm$.

In 2014 it was observed that certain discrete second-order Riesz transforms on products of integers are as well behaved as their continuous counterparts. In particular the discretised real part of the Beurling--Ahlfors operator,
\[
	B_{\mathfrak{r}}
	=R_1^2-R_2^2
	=R_1^+R_1^--R_2^+R_2^-,
\]
has the same optimal bound~$p^\star-1$ as in the continuous case (the upper bound of Nazarov--Volberg~\cite{NV2004}, with sharpness due to Geiss--Montgomery-Smith--Saksman~\cite{GMSS2009}). The proof in~\cite{DP2014} is deterministic and uses Bellman functions. Thanks to the specific form of the operator, the dualised form did not require mixed derivatives and permitted an estimate from the mere existence of a Bellman function, without an explicit formula. Later, Arcozzi--Domelevo--Petermichl~\cite{ArcDomPet2016a} treated second-order Riesz transforms in a probabilistic semi-discrete setting. That work generalised the deterministic results of~\cite{DP2014} and gave an ideological explanation of why the sharp constant~$p^\star-1$ was available: strong subordination holds and allows one to apply Burkholder's estimate directly, or rather its continuous-time analogue from Bañuelos--Wang~\cite{BW1995}. For this reason both~\cite{DP2014} and~\cite{ArcDomPet2016a} were limited to discrete second-order Riesz transforms of the type~$\sum_{i=1}^N\alpha_i R_i^2$ with~$|\alpha_i|\le 1$. In the deterministic picture, strong differential subordination translates into a straight use of the Bellman dynamics via a rough estimate of the Hessian along straight-line paths in the Taylor formula over a convex domain. That feature is specific to the second-order transforms considered in~\cite{DP2014}.

In 2018 Bañuelos and Kwaśnicki released their remarkable sharp estimate for the naïve discrete Hilbert transform~\cite{BK2019}. They realise that operator as the conditional expectation of a martingale transform built from a Doob~$h$-process in the upper half-plane. This produces an auxiliary convolution operator whose~$\ell^p$-norm equals~$\cot\bigl(\pi/(2p^\star)\bigr)$ and which differs from the discrete Hilbert transform by convolution against a probability kernel, whence the same sharp cotangent bound follows. Their auxiliary operator is free of jumps and feeds directly into the known estimates that use Essén's function. The strategy does not easily generalise to other operators or to higher dimensions.

Heat extensions in the upper half-space, originating in Petermichl--Volberg~\cite{PV2002} and used by Nazarov--Volberg~\cite{NV2004}, give~$L^p$ estimates for second-order Riesz transforms based on Burkholder's theorems, as part of the then-best estimate for the Beurling--Ahlfors operator, whose real and imaginary parts are themselves second-order Riesz transforms.

We reach our goal via a probabilistic model and an identity involving stochastic integration, in the spirit of~\cite{ArcDomPet2016a}. We use stochastic integrals with jump components adapted to functions on a product of integers with a deterministic time coordinate, and we show that the corresponding second-order Riesz transforms applied to a function may be written as the conditional expectation of a simple transformation of a stochastic integral associated with that function. Unlike~\cite{ArcDomPet2016a}, the martingale multiplier is not represented by a diagonal matrix. The desired~$\ell^p$ estimates therefore do not follow from Wang~\cite{Wang1995}, who identified the correct requirements on differential subordination in the presence of spatial discontinuities. Our pairs of martingales do not enjoy strong differential subordination. We reroute the jumps in such a way that the Burkholder functional can still be used, albeit indirectly.

Even in the continuous world, the conjectured estimate for the~$L^p$ norms of the Beurling--Ahlfors transform remains a fascinating open question. The argument of~\cite{NV2004} is deterministic in spirit and gave a new best estimate at the time, twice larger than expected. Subsequent improvements profited from the finer structure that survives in stochastic integration and from the probabilistic representation of second-order Riesz transforms due to Bañuelos--Méndez-Hernández~\cite{BMH2003}. That representation was used on many occasions, notably in the improvements of Bañuelos--Janakiraman~\cite{BJ2008} and Borichev--Janakiraman--Volberg~\cite{BJV2011}, through refined stochastic formulae and a new martingale estimate. The use of orthogonality in a jump process is at present out of reach, and the estimate~$2(p^\star-1)$ is the best currently available for the full discrete Beurling--Ahlfors operator.

The sharpness of the constant for the real part of the Beurling--Ahlfors operator---a combination of perfect squares of Riesz transforms---was proved by probabilistic methods in conjunction with a modification of a technique of Bourgain in Geiss--Montgomery-Smith--Saksman~\cite{GMSS2009}. The real part alone already attains the conjectured norm.

The important discretised Beurling--Ahlfors operator is excluded from those considerations. Its imaginary part on~$\mathbb{R}^2$ is~$-2R_1R_2$. In the continuous world this operator is obtained from the real part~$R_1^2-R_2^2$ by a rotation of~$\pi/4$, so the estimate follows immediately. On integer grids the rotation does not leave the grid invariant.

In this paper we derive a stochastic formula for the discrete imaginary part
\[
	B_{\mathfrak{i}}
	=-R_1^+R_2^--R_1^-R_2^+
\]
and prove that its sharp~$\ell^p$ norm is still~$p^\star-1$. Unlike the cases treated in~\cite{ArcDomPet2016a}, the pair of martingales representing~$B_{\mathfrak{i}}f$ and~$f$ does not have the strong subordination property required for a direct application of Burkholder's estimate. We circumvent the obstacle by rerouting jumps through fictitious intermediate points, in a manner that distantly resembles a discrete-time construction of Treil~\cite{Tre2013b}. In continuous time the argument is different and requires a careful treatment of the compensators of the Poisson jump process.

The estimate generalises to higher dimensions and to certain second-order Riesz transforms whose martingale representations lack differential subordination. The bound remains~$p^\star-1$, independently of the dimension.

On a related note we mention dimension-free (but not sharp) estimates for the first-order Riesz vector. Even in the continuous case the sharp estimate is unknown. Going back to 2004, Lust-Piquard~\cite{LP2004} proved a remarkable dimension-free estimate for the discrete Riesz vector on products of integers for~$p\ge 2$, using noncommutative analysis. That work was carried further by Junge, Mei and Parcet~\cite{JMP2018}. There is a known obstruction for~$1<p<2$, where the discrete Riesz vector does \emph{not} have dimension-free bounds. In the noncommutative theory the obstruction is visible in a noncommutative Khintchine inequality. In the continuous case a dimension-free estimate is available for the full range of~$p$, including a proof that uses subordinate martingales. Only recently, Ivanisvili, Volberg and the authors~\cite{DIPV2026} recovered the dimension-free estimate on the lattice by a Bellman-function argument. It had been unclear for a long time whether the noncommutative methods could be circumvented at all. In that proof the obstruction for~$p<2$ is read off the Nazarov--Treil Bellman function~\cite{NazTre1996a}, which changes with the range of~$p$. Substantial new ideas were needed to recover the estimate. Those ideas do not apply to the imaginary part of the discretised Beurling--Ahlfors operator considered here, and would in any case not give optimal bounds.

\section{Notations and main results}

\subsection{Notations and operators}
\label{sec:notation}

Let~$e_1,\ldots,e_N$ denote the canonical basis of~$\mathbb{Z}^N$, so that
\[
	(e_i)_j=\delta_{ij},\qquad i,j=1,\ldots,N.
\]
For a function~$f:\mathbb{Z}^N\to\mathbb{R}$ the forward and backward differences in the~$i$-th coordinate are
\begin{align}
	\partial_i^+f(x)
	 & =f(x+e_i)-f(x),
	\label{eq:forward} \\
	\partial_i^-f(x)
	 & =f(x)-f(x-e_i).
	\label{eq:backward}
\end{align}
These two first-order operators commute, and the centred second difference is independent of order. We write
\begin{equation}
	\label{eq:second-diff}
	\begin{aligned}
		\partial_i^2f(x)
		 & :=\partial_i^+\partial_i^-f(x)
		=\partial_i^-\partial_i^+f(x)
		=f(x+e_i)-2f(x)+f(x-e_i)               \\
		 & =\partial_i^+f(x)-\partial_i^-f(x),
	\end{aligned}
\end{equation}
and the companion sum of one-sided derivatives
\[
	\partial_i^0f(x)
	:=\partial_i^+f(x)+\partial_i^-f(x).
\]
The discrete Laplacian (graph Laplacian of the integer grid) is
\begin{equation}
	\label{eq:laplacian}
	\Delta f(x)
	=\sum_{i=1}^N\partial_i^2f(x)
	=\sum_{y\sim x}\bigl(f(y)-f(x)\bigr),
\end{equation}
where~$y\sim x$ means that~$y$ is a nearest neighbour of~$x$. As an operator on~$\ell^2(\mathbb{Z}^N)$ one has~$\Delta\le 0$, with Fourier multiplier
\[
	\widehat{\Delta}(\theta)
	=\sum_{i=1}^N\bigl(2\cos\theta_i-2\bigr)
	=-4\sum_{i=1}^N\sin^2(\theta_i/2),\qquad\theta\in\mathbb{T}^N.
\]
The constant functions are formally harmonic, but they do not belong to~$\ell^2(\mathbb{Z}^N)$; in particular~$0$ is not an eigenvalue of~$\Delta$ on~$\ell^2$.

The discrete Riesz transforms associated with the two one-sided gradients are
\begin{equation}
	\label{eq:riesz}
	R_i^+
	=\partial_i^+\circ(-\Delta)^{-1/2},
	\qquad
	R_i^-
	=\partial_i^-\circ(-\Delta)^{-1/2},
\end{equation}
understood via the functional calculus of~$-\Delta$ (or equivalently as Fourier multipliers). They differ from one another by a unimodular trigonometric factor and are bounded on~$\ell^p(\mathbb{Z}^N)$ for~$1<p<\infty$. All of these operators commute, whatever the coordinate we consider. Namely
\begin{equation}
	\label{eq:commutation}
	\forall\,\alpha,\beta\in\{+,-\},\;
	\forall\,i,j\in\{1,\ldots,N\},\qquad
	\partial_i^\alpha\partial_j^\beta=\partial_j^\beta\partial_i^\alpha,
	\qquad
	R_i^\alpha R_j^\beta=R_j^\beta R_i^\alpha.
\end{equation}
In particular the diagonal second-order transforms may be written
\begin{equation}
	\label{eq:riesz2diag}
	R_i^2
	:=R_i^+R_i^-
	=R_i^-R_i^+
	=\frac12\bigl(R_i^+R_i^-+R_i^-R_i^+\bigr)
	=\partial_i^2\circ(-\Delta)^{-1},
\end{equation}
and similarly the imaginary part of the Beurling--Ahlfors operator can be written in its symmetrised form
\begin{equation}
	\label{eq:Bi-sym}
	B_{\mathfrak{i}}
	=-\frac12\bigl(R_1^+R_2^-+R_2^+R_1^-+R_1^-R_2^++R_2^-R_1^+\bigr).
\end{equation}

The~$\ell^2$ pairing on the lattice is written
\[
	\langle f,g\rangle
	=\sum_{x\in\mathbb{Z}^N}f(x)\,g(x).
\]
The \emph{augmented discrete gradient} collects every one-sided difference at a point. We order the components by writing first all forward differences and then all backward differences,
\begin{equation}
	\label{eq:auggrad}
	\widetilde{\nabla}f(x)
	=\bigl(
	\partial_1^+f(x),\,
	\ldots,\,
	\partial_N^+f(x),\,
	\partial_1^-f(x),\,
	\ldots,\,
	\partial_N^-f(x)
	\bigr)^\top
	\in\mathbb{R}^{2N}.
\end{equation}
It takes values in an augmented discrete tangent space of dimension~$2N$, twice the combinatorial dimension of the grid. The associated pairing is
\begin{equation}
	\label{eq:augpair}
	\langle\widetilde{\nabla}f,\widetilde{\nabla}g\rangle
	=\sum_{x\in\mathbb{Z}^N}
	\sum_{i=1}^N
	\Bigl(
	\partial_i^+f(x)\,\partial_i^+g(x)
	+\partial_i^-f(x)\,\partial_i^-g(x)
	\Bigr).
\end{equation}
If~$\mathbf{e}_1,\ldots,\mathbf{e}_{2N}$ denotes the canonical basis of~$\mathbb{R}^{2N}$, one has the identifications
\[
	\partial_i^+f=\mathbf{e}_i^\top\widetilde{\nabla}f,
	\qquad
	\partial_i^-f=\mathbf{e}_{N+i}^\top\widetilde{\nabla}f,
	\qquad
	i=1,\ldots,N.
\]
Summation by parts gives the discrete Green identities
\begin{equation}
	\label{eq:green}
	-\langle\Delta f,g\rangle
	=\sum_{i=1}^N\langle\partial_i^+f,\partial_i^+g\rangle
	=\sum_{i=1}^N\langle\partial_i^-f,\partial_i^-g\rangle
	=\frac12\langle\widetilde{\nabla}f,\widetilde{\nabla}g\rangle,
\end{equation}
valid whenever~$f,g$ and the first differences are square-summable (e.g.\ for finitely supported functions, and then by density on the form domain of~$\Delta$).

\subsection{Second order Riesz transforms and main results}
\label{sec:second}

When~$N=2$, the real and imaginary parts of the discrete Beurling--Ahlfors operator on~$\mathbb{Z}^2$ are
\begin{equation}
	\label{eq:beurling}
	\begin{aligned}
		B_{\mathfrak{r}}
		 & =R_1^2-R_2^2
		=\frac12\bigl(R_1^+R_1^-+R_1^-R_1^+\bigr)
		-\frac12\bigl(R_2^+R_2^-+R_2^-R_2^+\bigr),                           \\
		B_{\mathfrak{i}}
		 & =-\frac12\bigl(R_1^+R_2^-+R_2^+R_1^-+R_1^-R_2^++R_2^-R_1^+\bigr).
	\end{aligned}
\end{equation}
Among others, we prove the exact norm estimate for the imaginary part. The sharp bound for the real part was already proved in Arcozzi--Domelevo--Petermichl~\cite{ArcDomPet2016a}, and will be a consequence of Theorem~\ref{thm:main} below. Write
\[
	p^\star-1
	:=\max\bigl\{p-1,(p-1)^{-1}\bigr\},
	\qquad
	1<p<\infty.
\]
Our first goal is to prove:

\begin{theorem}
	\label{thm:imag}
	For~$1<p<\infty$ one has the sharp bound
	\[
		\|B_{\mathfrak{i}}\|_{\ell^p(\mathbb{Z}^2)\to\ell^p(\mathbb{Z}^2)}
		=p^\star-1.
	\]
\end{theorem}

The continuous analogue of this estimate on~$\mathbb{R}^2$ is due to Nazarov and Volberg~\cite{NV2004}; in the continuous case only, it follows from the corresponding estimate for the diagonal second-order Riesz transform~$B_{\mathfrak{r}}$ by observing that~$B_{\mathfrak{i}}$ is obtained from~$B_{\mathfrak{r}}$ through a rotation by~$\pi/4$. In the discrete case the nature of these estimates is distinct: the rotation argument is unavailable, because the integer grid is not invariant under rotations by angles other than multiples of~$\pi/2$.

\begin{remark}
	\label{rem:sharp-imag}
	The sharpness in Theorem~\ref{thm:imag} is a consequence of the sharpness for~$R_1^2-R_2^2$ proved by Geiss--Montgomery-Smith--Saksman~\cite{GMSS2009}. As mentioned in the introduction, the continuous imaginary part is obtained from the real part~$R_1^2-R_2^2$ by a rotation of~$\pi/4$, so the sharp estimate follows immediately for the continuous operator~$B_{\mathfrak{i}}$. The sharpness for the discrete version is obtained in the standard way by mesh refinement: the continuous operator arises as a scaling limit of the discrete one, and therefore cannot have a strictly smaller norm. This comparison of discrete and continuous norms by passing to the continuum is the same procedure recorded for the Hilbert transform by Laeng~\cite{Laeng2007} and recalled in the introduction of Bañuelos--Kwaśnicki~\cite{BK2019}.
\end{remark}

Our arguments allow for a more general estimate that includes the estimate for~$B_{\mathfrak{i}}$. Let~$\mathcal{A}$ be any subset of~$\{1,\ldots,N\}^2$ such that
\begin{itemize}[leftmargin=*,itemsep=2pt]
	\item $(i,j)\in\mathcal{A}$ implies~$(j,i)\in\mathcal{A}$,
	\item $(i,j)\in\mathcal{A}$ implies~$i\neq j$,
	\item $(i,j)\in\mathcal{A}$ implies~$(i,j')\notin\mathcal{A}$ for every~$j'\neq j$.
\end{itemize}
In other words,~$\mathcal{A}$ is a symmetric set of off-diagonal ordered pairs in which each index appears in at most one pair: for each~$i$ there is at most one~$j$ with~$(i,j)\in\mathcal{A}$. We will say~$\{i,j\}\in\mathcal{A}$ if~$(i,j)\in\mathcal{A}$ or~$(j,i)\in\mathcal{A}$. Write
\[
	\mathcal{A}'
	=\bigl\{i\in\{1,\ldots,N\}:(i,j)\notin\mathcal{A}\text{ for all }j\bigr\}
\]
for the complementary set of unmatched indices. For~$i\neq j$ define the fully symmetrised off-diagonal transform
\begin{equation}
	\label{eq:offdiag}
	\begin{aligned}
		R_{(i,j)}^2
		 & :=\frac14\bigl(R_i^+R_j^-+R_j^+R_i^-+R_i^-R_j^++R_j^-R_i^+\bigr) \\
		 & =\frac12\bigl(R_i^+R_j^-+R_j^+R_i^-\bigr)
		=\frac12\bigl(R_i^-R_j^++R_j^-R_i^+\bigr).
	\end{aligned}
\end{equation}
(The two ordered pairs~$(i,j)$ and~$(j,i)$ determine the same operator.) The same formula with~$j=i$ recovers the diagonal transforms, and the imaginary part of the Beurling--Ahlfors operator is a multiple of the unique off-diagonal pair on~$\mathbb{Z}^2$:
\begin{equation}
	\label{eq:diag-from-sym}
	R_i^2=R_{(i,i)}^2,
	\qquad
	B_{\mathfrak{i}}=-2R_{(1,2)}^2.
\end{equation}
Consider, for any such set~$\mathcal{A}$ and associated~$\mathcal{A}'$, and for any sequence~$(\alpha_i)_{1\le i\le N}$ with~$|\alpha_i|\le 1$ and with~$\alpha_i=\alpha_j$ whenever~$(i,j)\in\mathcal{A}$, the second-order Riesz transform
\begin{equation}
	\label{eq:RalphaA}
	R_{\alpha,\mathcal{A}}^2
	=\sum_{i\in\mathcal{A}'}\alpha_i\,R_{(i,i)}^2
	+\sum_{\{i,j\}\in\mathcal{A}}2\alpha_i\,R_{(i,j)}^2.
\end{equation}
This operator admits a compact writing in terms of the augmented discrete Riesz vector. In the same ordering as~\eqref{eq:auggrad} we set
\begin{equation}
	\label{eq:augriesz}
	\widetilde{R}
	:=\bigl(R_1^+,\ldots,R_N^+,R_1^-,\ldots,R_N^-\bigr)^\top
	=\widetilde{\nabla}\circ(-\Delta)^{-1/2}.
\end{equation}
The adjoint of a one-sided Riesz transform exchanges the forward and backward differences and changes the sign:
\[
	\bigl(R_i^+\bigr)^\ast=-R_i^-,
	\qquad
	\bigl(R_i^-\bigr)^\ast=-R_i^+,
\]
and likewise~$\bigl(\partial_i^+\bigr)^\ast=-\partial_i^-$,~$\bigl(\partial_i^-\bigr)^\ast=-\partial_i^+$. Consequently the adjoint of the augmented vector~$\widetilde{R}$ is the row obtained by swapping the two blocks of length~$N$ and inserting a minus. Each first-order component is a coordinate of the augmented gradient:
\[
	\partial_i^+=\mathbf{e}_i^\top\widetilde{\nabla},
	\qquad
	\partial_i^-=\mathbf{e}_{N+i}^\top\widetilde{\nabla},
	\qquad
	R_i^+=\mathbf{e}_i^\top\widetilde{R},
	\qquad
	R_i^-=\mathbf{e}_{N+i}^\top\widetilde{R}.
\]
The elementary building blocks of~\eqref{eq:offdiag} are combinations of operators of the form~$D^\ast D'$, where~$D,D'$ are first-order Riesz transforms. In the diagonal case this is a genuine square: \(R_i^2=-(R_i^+)^\ast R_i^+=-(R_i^-)^\ast R_i^-\). In the off-diagonal case one has the polarisation
\[
	R_{(i,j)}^2
	=-\frac14\sum\bigl((R_i^\pm)^\ast R_j^\pm+(R_j^\pm)^\ast R_i^\pm\bigr).
\]
Assembling these identities, the operator~\eqref{eq:RalphaA} is the quadratic form
\begin{equation}
	\label{eq:R-M-R}
	R_{\alpha,\mathcal{A}}^2
	=\widetilde{R}^*\,\mathbb{M}_{\alpha,\mathcal{A}}\,\widetilde{R},
\end{equation}
where~$\mathbb{M}_{\alpha,\mathcal{A}}$ is the symmetric~$2N\times 2N$ matrix
\begin{equation}
	\label{eq:Malpha}
	\begin{aligned}
		\mathbb{M}_{\alpha,\mathcal{A}}
		 & =
		-\sum_{i\in\mathcal{A}'}\frac{\alpha_i}{2}
		\bigl(\mathbf{e}_i\mathbf{e}_i^\ast+\mathbf{e}_{N+i}\mathbf{e}_{N+i}^\ast\bigr) \\
		 & \qquad
		-\sum_{\{i,j\}\in\mathcal{A}}\frac{\alpha_i}{2}
		\bigl(
		\mathbf{e}_i\mathbf{e}_j^\ast+\mathbf{e}_j\mathbf{e}_i^\ast
		+\mathbf{e}_{N+i}\mathbf{e}_{N+j}^\ast+\mathbf{e}_{N+j}\mathbf{e}_{N+i}^\ast
		\bigr).
	\end{aligned}
\end{equation}
Thus~$\mathbb{M}_{\alpha,\mathcal{A}}$ is a linear combination of the coefficients~$\alpha_i$ with the elementary matrices~$\mathbf{e}_k\mathbf{e}_\ell^\ast$. It acts separately and identically on the forward block and on the backward block. In particular~$|\mathbb{M}_{\alpha,\mathcal{A}}|_{\mathrm{op}}\le\tfrac12$ whenever~$|\alpha|_\infty\le 1$.

The two model operators on~$\mathbb{Z}^2$ are recovered as
\[
	B_{\mathfrak{r}}
	=R_{\alpha,\emptyset}^2
	\quad\text{with}\quad\alpha=(1,-1)
\]
and
\[
	B_{\mathfrak{i}}
	=-2R_{(1,2)}^2
	=R_{\alpha,\mathcal{A}}^2
	\quad\text{with}\quad
	\mathcal{A}=\bigl\{(1,2),(2,1)\bigr\}
	\text{ and }
	\alpha=(-1,-1).
\]

\begin{theorem}
	\label{thm:main}
	For any set~$\mathcal{A}$ as above and any coefficient sequence~$\alpha$ as above with~$|\alpha|_\infty\le 1$, one has the sharp bound
	\[
		\bigl\|R_{\alpha,\mathcal{A}}^2\bigr\|_{\ell^p(\mathbb{Z}^N)\to\ell^p(\mathbb{Z}^N)}
		\le p^\star-1.
	\]
\end{theorem}

Theorem~\ref{thm:imag} is the special case~$N=2$,~$\mathcal{A}=\{(1,2),(2,1)\}$,~$\alpha=(-1,-1)$, together with the matching lower bound inherited from the continuous operator by mesh refinement. Though Theorem~\ref{thm:imag} is a corollary of Theorem~\ref{thm:main}, we will rather prove first Theorem~\ref{thm:imag} and observe how the method of proof extends to the operators considered in Theorem~\ref{thm:main}.

\section{Heat extensions and weak formulations}
\label{sec:heat}

\subsection{The reproducing formula}

The semi-discrete heat equation on the lattice is the Cauchy problem
\begin{equation}
	\label{eq:heat}
	\partial_t\widetilde{f}(t,x)=\Delta\widetilde{f}(t,x),\qquad
	\widetilde{f}(0,x)=f(x),\qquad
	x\in\mathbb{Z}^N,\ t>0.
\end{equation}
Its solution is the heat extension
\[
	\widetilde{f}(t,x)=(e^{t\Delta}f)(x).
\]
The same construction with initial datum~$g$ is written~$\widetilde{g}=e^{t\Delta}g$. Because~$\Delta\le 0$ on~$\ell^2$ and~$0$ is not an eigenvalue, one has~$\widetilde{f}(t,\,\cdot\,)\to 0$ in~$\ell^2(\mathbb{Z}^N)$ as~$t\to\infty$, for every~$f\in\ell^2(\mathbb{Z}^N)$.

Differentiating the pairing in time and applying~\eqref{eq:green} yields
\begin{align*}
	\frac{d}{dt}\langle\widetilde{f}(t,\,\cdot\,),\widetilde{g}(t,\,\cdot\,)\rangle
	 & =2\langle\Delta\widetilde{f}(t,\,\cdot\,),\widetilde{g}(t,\,\cdot\,)\rangle                                \\
	 & =-2\sum_{i=1}^N\langle\partial_i^+\widetilde{f}(t,\,\cdot\,),\partial_i^+\widetilde{g}(t,\,\cdot\,)\rangle \\
	 & =-2\sum_{i=1}^N\langle\partial_i^-\widetilde{f}(t,\,\cdot\,),\partial_i^-\widetilde{g}(t,\,\cdot\,)\rangle \\
	 & =-\langle\widetilde{\nabla}\widetilde{f}(t,\,\cdot\,),\widetilde{\nabla}\widetilde{g}(t,\,\cdot\,)\rangle.
\end{align*}
Integrating from~$0$ to~$\infty$ and using the decay at infinity gives the weak reproducing formula, written in the equivalent forms
\begin{equation}
	\label{eq:reproducing}
	\begin{aligned}
		\langle f,g\rangle
		 & =2\int_0^\infty\sum_{i=1}^N\sum_{x\in\mathbb{Z}^N}
		\partial_i^+\widetilde{f}(t,x)\,\partial_i^+\widetilde{g}(t,x)\,dt \\
		 & =2\int_0^\infty\sum_{x\in\mathbb{Z}^N}\sum_{i=1}^N
		\partial_i^+\widetilde{f}(t,x)\,\partial_i^+\widetilde{g}(t,x)\,dt \\
		 & =\int_0^\infty\sum_{x\in\mathbb{Z}^N}\sum_{i=1}^N
		\Bigl(
		\partial_i^+\widetilde{f}(t,x)\,\partial_i^+\widetilde{g}(t,x)
		+\partial_i^-\widetilde{f}(t,x)\,\partial_i^-\widetilde{g}(t,x)
		\Bigr)\,dt                                                         \\
		 & =\int_0^\infty
		\langle\widetilde{\nabla}\widetilde{f}(t,\,\cdot\,),\widetilde{\nabla}\widetilde{g}(t,\,\cdot\,)\rangle\,dt.
	\end{aligned}
\end{equation}
The first two lines may equally well be written with every~$\partial_i^+$ replaced by~$\partial_i^-$. Equivalently, in operator form,
\begin{equation}
	\label{eq:reproducing-lap}
	\langle f,g\rangle
	=-2\int_0^\infty\langle\Delta\widetilde{f}(t,\,\cdot\,),\widetilde{g}(t,\,\cdot\,)\rangle\,dt.
\end{equation}
All of these formulae extend from finitely supported data to all of~$\ell^2(\mathbb{Z}^N)$ by density on the form domain of~$\Delta$.

\subsection{Weak formulations for second-order Riesz transforms}
\label{sec:weak-second}

We derive the elementary pairings from~\eqref{eq:reproducing-lap}, the definitions~\eqref{eq:riesz} and the commutation relations~\eqref{eq:commutation}. The heat semigroup commutes with~$\Delta$ and with every~$\partial_i^\pm$ and~$R_i^\pm$, so~$\widetilde{Tf}=T\widetilde{f}$ for each of these operators. Applying~\eqref{eq:reproducing-lap} to the pair~$\bigl(R_i^+R_j^-f,\,g\bigr)$ therefore gives
\[
	\langle R_i^+R_j^-f,g\rangle
	=-2\int_0^\infty\bigl\langle\Delta R_i^+R_j^-\widetilde{f}(t,\,\cdot\,),\,\widetilde{g}(t,\,\cdot\,)\bigr\rangle\,dt.
\]
By definition and commutation,~$R_i^+R_j^-=\partial_i^+\partial_j^-\circ(-\Delta)^{-1}$. Since~$\Delta\circ(-\Delta)^{-1}=-\mathrm{Id}$ on~$\ell^2(\mathbb{Z}^N)$,
\[
	\Delta R_i^+R_j^-\widetilde{f}
	=-\partial_i^+\partial_j^-\widetilde{f}.
\]
The adjoint identity~$\bigl(\partial_i^+\bigr)^\ast=-\partial_i^-$ then yields
\[
	\bigl\langle\partial_i^+\partial_j^-\widetilde{f},\,\widetilde{g}\bigr\rangle
	=\bigl\langle\partial_j^-\widetilde{f},\,\bigl(\partial_i^+\bigr)^\ast\widetilde{g}\bigr\rangle
	=-\bigl\langle\partial_j^-\widetilde{f},\,\partial_i^-\widetilde{g}\bigr\rangle,
\]
and we obtain
\begin{equation}
	\label{eq:weak-pm}
	\begin{aligned}
		\langle R_i^+R_j^-f,g\rangle
		 & =-2\int_0^\infty\sum_{x\in\mathbb{Z}^N}
		\partial_j^-\widetilde{f}(t,x)\,\partial_i^-\widetilde{g}(t,x)\,dt, \\
		\langle R_i^-R_j^+f,g\rangle
		 & =-2\int_0^\infty\sum_{x\in\mathbb{Z}^N}
		\partial_j^+\widetilde{f}(t,x)\,\partial_i^+\widetilde{g}(t,x)\,dt.
	\end{aligned}
\end{equation}
The second line is the same computation with the roles of the forward and backward differences exchanged: now~$\bigl(R_i^-\bigr)^\ast=-R_i^+$, so the pairing becomes a product of two forward transforms.

The off-diagonal operators~\eqref{eq:offdiag} are averages of these elementary products. Using the backward representation~$R_{(i,j)}^2=\tfrac12\bigl(R_i^+R_j^-+R_j^+R_i^-\bigr)$ and the first line of~\eqref{eq:weak-pm} (together with the same line after~$i\leftrightarrow j$) gives the first identity below; the forward representation of~\eqref{eq:offdiag} gives the second; averaging the two yields the third:
\begin{equation}
	\label{eq:weak-off}
	\begin{aligned}
		\langle R_{(i,j)}^2f,g\rangle
		 & =-\int_0^\infty\sum_{x\in\mathbb{Z}^N}
		\Bigl(
		\partial_i^-\widetilde{f}\,\partial_j^-\widetilde{g}
		+\partial_j^-\widetilde{f}\,\partial_i^-\widetilde{g}
		\Bigr)\,dt                                       \\
		 & =-\int_0^\infty\sum_{x\in\mathbb{Z}^N}
		\Bigl(
		\partial_i^+\widetilde{f}\,\partial_j^+\widetilde{g}
		+\partial_j^+\widetilde{f}\,\partial_i^+\widetilde{g}
		\Bigr)\,dt                                       \\
		 & =-\frac12\int_0^\infty\sum_{x\in\mathbb{Z}^N}
		\Bigl(
		\partial_i^+\widetilde{f}\,\partial_j^+\widetilde{g}
		+\partial_j^+\widetilde{f}\,\partial_i^+\widetilde{g}
		+\partial_i^-\widetilde{f}\,\partial_j^-\widetilde{g}
		+\partial_j^-\widetilde{f}\,\partial_i^-\widetilde{g}
		\Bigr)\,dt.
	\end{aligned}
\end{equation}
The diagonal operators of~\eqref{eq:RalphaA} are the special case~$j=i$. With the notation~$R_{(i,i)}^2=R_i^2$ one finds
\begin{equation}
	\label{eq:weak-diag}
	\begin{aligned}
		\langle R_{(i,i)}^2f,g\rangle
		 & =-2\int_0^\infty\sum_{x\in\mathbb{Z}^N}
		\partial_i^+\widetilde{f}(t,x)\,\partial_i^+\widetilde{g}(t,x)\,dt \\
		 & =-2\int_0^\infty\sum_{x\in\mathbb{Z}^N}
		\partial_i^-\widetilde{f}(t,x)\,\partial_i^-\widetilde{g}(t,x)\,dt \\
		 & =-\int_0^\infty\sum_{x\in\mathbb{Z}^N}
		\Bigl(
		\partial_i^+\widetilde{f}(t,x)\,\partial_i^+\widetilde{g}(t,x)
		+\partial_i^-\widetilde{f}(t,x)\,\partial_i^-\widetilde{g}(t,x)
		\Bigr)\,dt.
	\end{aligned}
\end{equation}
Collecting~\eqref{eq:weak-off} and~\eqref{eq:weak-diag} into a single pairing against the augmented gradient yields the following statement.

\begin{lemma}
	\label{lem:weak-general}
	Let~$\mathbb{M}_{\alpha,\mathcal{A}}$ be the matrix~\eqref{eq:Malpha}. Then
	\begin{equation}
		\label{eq:weak-general}
		\langle R_{\alpha,\mathcal{A}}^2f,g\rangle
		=2\int_0^\infty
		\sum_{x\in\mathbb{Z}^N}
		\bigl\langle
		\mathbb{M}_{\alpha,\mathcal{A}}\,
		\widetilde{\nabla}\widetilde{f}(t,x),\,
		\widetilde{\nabla}\widetilde{g}(t,x)
		\bigr\rangle
		\,dt.
	\end{equation}
	In particular~$|\mathbb{M}_{\alpha,\mathcal{A}}|_{\mathrm{op}}\le\tfrac12$ whenever~$|\alpha|_\infty\le 1$.
\end{lemma}

The factor~$2$ converts the coefficients~$\pm\alpha_i/2$ of~\eqref{eq:Malpha} into the heat pairings~\eqref{eq:weak-diag} and~\eqref{eq:weak-off}: for an unmatched index one has~$2\cdot(-\alpha_i/2)=-\alpha_i$, and for a pair~$\{i,j\}$ the same factor turns~$2\alpha_i R_{(i,j)}^2$ into the four-term identity of~\eqref{eq:weak-off}.

The imaginary part of the discrete Beurling--Ahlfors operator is the model case~$N=2$,~$\mathcal{A}=\{(1,2),(2,1)\}$,~$\alpha=(-1,-1)$, so that~$B_{\mathfrak{i}}=-2R_{(1,2)}^2$. With the ordering~\eqref{eq:auggrad} one has~$\partial_1^+=\mathbf{e}_1^\top\widetilde{\nabla}$,~$\partial_2^+=\mathbf{e}_2^\top\widetilde{\nabla}$,~$\partial_1^-=\mathbf{e}_3^\top\widetilde{\nabla}$,~$\partial_2^-=\mathbf{e}_4^\top\widetilde{\nabla}$, and Lemma~\ref{lem:weak-general} becomes the following statement.

\begin{lemma}
	\label{lem:R12}
	Write the augmented discrete gradient on~$\mathbb{Z}^2$ as
	\[
		\widetilde{\nabla}
		=\bigl(\partial_1^+,\partial_2^+,\partial_1^-,\partial_2^-\bigr)^\top.
	\]
	Then
	\[
		\langle B_{\mathfrak{i}}f,g\rangle
		=\sum_{x\in\mathbb{Z}^2}
		\int_0^\infty
		\left\langle
		\begin{pmatrix}
			0 & 1 & 0 & 0 \\
			1 & 0 & 0 & 0 \\
			0 & 0 & 0 & 1 \\
			0 & 0 & 1 & 0
		\end{pmatrix}
		\widetilde{\nabla}\widetilde{f}(t,x),\,
		\widetilde{\nabla}\widetilde{g}(t,x)
		\right\rangle
		\,dt.
	\]
\end{lemma}

The swap matrix in the lemma is~$2\mathbb{M}_{\alpha,\mathcal{A}}$ for~$\mathcal{A}=\{(1,2),(2,1)\}$ and~$\alpha=(-1,-1)$.

\section{Stochastic representations}
\label{sec:stoch}

\subsection{Compound Poisson process and Itô formula}
\label{sec:poisson}

The starting point is the weak formulation of Lemma~\ref{lem:R12}. We now put a compound Poisson process under this identity, following Arcozzi--Domelevo--Petermichl~\cite{ArcDomPet2016a}, on the full lattice~$\mathbb{Z}^N$. Let~$\mathcal{N}^1,\ldots,\mathcal{N}^N$ be independent Poisson processes of intensity~$2$, and write
\[
	\widetilde{\mathcal{N}}
	=\bigl(\mathcal{N}^1,\ldots,\mathcal{N}^N\bigr)
\]
for the vector process. The scalar process
\[
	\mathcal{N}_t
	:=\sum_{i=1}^N\mathcal{N}^i_t
\]
counts all jumps of~$\widetilde{\mathcal{N}}$ up to time~$t$. Let~$(\tau_k)_{k\ge 1}$ be i.i.d.\ signs with~$\mathbb{P}(\tau=\pm 1)=\tfrac12$, independent of~$\widetilde{\mathcal{N}}$; the sign attached to a jump at time~$s$ is~$\tau_{\mathcal{N}_s}$, indexed by the total jump counter. The lattice-valued process~$\widetilde{X}=(X^1,\ldots,X^N)$ on~$\mathbb{Z}^N$ is the compound Poisson process
\[
	\mathrm{d}X^i_s
	=\tau_{\mathcal{N}_s}\,\mathrm{d}\mathcal{N}^i_s,
	\qquad i=1,\ldots,N,
\]
so that~$\widetilde{X}$ jumps by~$\pm e_i$ at the jump times of~$\mathcal{N}^i$. The infinitesimal generator of~$\widetilde{X}$ is the discrete Laplacian~$\Delta$: each coordinate has total jump rate~$2$. In the direction~$e_i$ we use the difference and the sum of the one-sided derivatives already recorded in~\eqref{eq:second-diff},
\[
	\partial_i^2\varphi
	=\partial_i^+\varphi-\partial_i^-\varphi,
	\qquad
	\partial_i^0\varphi
	=\partial_i^+\varphi+\partial_i^-\varphi.
\]
The associated first-order fields on the tangent space are the~$\mathbb{R}^N$-valued vectors
\begin{equation}
	\label{eq:grad02}
	\widetilde{\nabla}^2\varphi
	=\bigl(\partial_1^2\varphi,\ldots,\partial_N^2\varphi\bigr),
	\qquad
	\widetilde{\nabla}^0\varphi
	=\bigl(\partial_1^0\varphi,\ldots,\partial_N^0\varphi\bigr).
\end{equation}
A jump of~$X^i$ with sign~$\tau$ changes a test function by
\[
	\varphi(x+\tau e_i)-\varphi(x)
	=\frac12\widetilde{\nabla}^2_i\varphi(x)+\frac\tau2\widetilde{\nabla}^0_i\varphi(x).
\]
Consequently, for a regular function~$\varphi=\varphi(t,x)$ on~$[0,\infty)\times\mathbb{Z}^N$, the Itô formula reads
\begin{equation}
	\label{eq:ito}
	\begin{aligned}
		\varphi(t,\widetilde{X}_t)-\varphi(0,\widetilde{X}_0)
		 & =\int_0^t
		\bigl(\partial_t\varphi(s,\widetilde{X}_{s-})+\Delta\varphi(s,\widetilde{X}_{s-})\bigr)\,\mathrm{d}s \\
		 & \qquad
		+\frac12\int_0^t
		\bigl\langle\widetilde{\nabla}^2\varphi(s,\widetilde{X}_{s-}),\,\mathrm{d}\bigl(\widetilde{\mathcal{N}}_s-2s\mathbf{1}\bigr)\bigr\rangle
		+\frac12\int_0^t
		\bigl\langle\widetilde{\nabla}^0\varphi(s,\widetilde{X}_{s-}),\,\mathrm{d}\widetilde{X}_s\bigr\rangle.
	\end{aligned}
\end{equation}
Here~$\mathbf{1}=(1,\ldots,1)\in\mathbb{R}^N$, so that~$\widetilde{\mathcal{N}}_s-2s\mathbf{1}$ is the vector of compensated Poisson processes, and the two scalar products are the Euclidean pairings on~$\mathbb{R}^N$. The last two integrals are martingales. In particular, if~$\varphi(t,x)=\widetilde{h}(T-t,x)$ is a backward heat extension, the finite-variation integrand vanishes and~$\varphi(t,\widetilde{X}_t)$ is a martingale.

The quadratic variations of the driving processes are
\begin{equation}
	\label{eq:qv-drivers}
	\mathrm{d}\bigl[\mathcal{N}^i-2s,\,\mathcal{N}^i-2s\bigr]_t
	=\mathrm{d}\mathcal{N}^i_t,
	\qquad
	\mathrm{d}\bigl[\mathcal{N}^i-2s,\,X^i\bigr]_t
	=\tau_{\mathcal{N}_t}\,\mathrm{d}\mathcal{N}^i_t,
	\qquad
	\mathrm{d}\bigl[X^i,X^i\bigr]_t
	=\mathrm{d}\mathcal{N}^i_t,
\end{equation}
and distinct coordinates do not jump together.

\subsection{Product of martingales and the Gundy--Varopoulos projection}
\label{sec:product}

We are now going to study martingale transforms corresponding to the matrices~$\mathbb{M}_{\alpha,\mathcal{A}}$ defined in~\eqref{eq:Malpha}. Following the Gundy--Varopoulos approach, the product of such a transform with a test martingale delivers the weak formulations of Lemma~\ref{lem:weak-general} and allows us to realise the second-order Riesz transforms under consideration as projection operators through conditional expectation with respect to arrival points.

For that we first study the effect of the elementary martingale transform associated with the elementary block~$\mathbb{M}_{i,j}$ appearing as a summand in the second sum of~\eqref{eq:Malpha}. For~$i\neq j$ we set
\begin{equation}
	\label{eq:Mij}
	\mathbb{M}_{i,j}
	:=-\frac12\bigl(
	\mathbf{e}_i\mathbf{e}_j^\ast+\mathbf{e}_j\mathbf{e}_i^\ast
	+\mathbf{e}_{N+i}\mathbf{e}_{N+j}^\ast+\mathbf{e}_{N+j}\mathbf{e}_{N+i}^\ast
	\bigr),
\end{equation}
so that the off-diagonal part of~\eqref{eq:Malpha} is~$\sum_{\{i,j\}\in\mathcal{A}}\alpha_i\,\mathbb{M}_{i,j}$.

Fix a large height~$T>0$. As in Gundy--Varopoulos~\cite{GV1976}, one starts an independent copy of the compound Poisson process from every lattice point: for each~$x\in\mathbb{Z}^N$ the background noise is the space-time path~$(T-t,\widetilde{X}_t)_{0\le t\le T}$ with~$\widetilde{X}_0=x$. Applying~\eqref{eq:ito} to the backward heat extensions of~$f$ and~$g$ produces the martingales. The test martingale is the heat extension itself,
\begin{align}
	\label{eq:Mg}
	M^g_t
	 & =\widetilde{g}(T-t,\widetilde{X}_t)
	=\widetilde{g}(T,\widetilde{X}_0)
	+\frac12\int_0^t
	\Bigl(
	\bigl\langle\widetilde{\nabla}^2\widetilde{g}(T-s,\widetilde{X}_{s-}),\,\mathrm{d}\bigl(\widetilde{\mathcal{N}}_s-2s\mathbf{1}\bigr)\bigr\rangle
	+\bigl\langle\widetilde{\nabla}^0\widetilde{g}(T-s,\widetilde{X}_{s-}),\,\mathrm{d}\widetilde{X}_s\bigr\rangle
	\Bigr),
\end{align}
and likewise~$M^f_t=\widetilde{f}(T-t,\widetilde{X}_t)$. The martingale transforms start at zero: we set
\begin{align}
	\label{eq:MAf}
	M^{\{i,j\},f}_t
	 & =\frac12\int_0^t
	\Bigl(
	\bigl\langle A_{i,j}\widetilde{\nabla}^2\widetilde{f}(T-s,\widetilde{X}_{s-}),\,\mathrm{d}\bigl(\widetilde{\mathcal{N}}_s-2s\mathbf{1}\bigr)\bigr\rangle
	+\bigl\langle A_{i,j}\widetilde{\nabla}^0\widetilde{f}(T-s,\widetilde{X}_{s-}),\,\mathrm{d}\widetilde{X}_s\bigr\rangle
	\Bigr),
\end{align}
where~$A_{i,j}$ is the~$N\times N$ matrix that swaps the coordinates~$i$ and~$j$ and vanishes on the complementary subspace,
\begin{equation}
	\label{eq:Aij}
	A_{i,j}
	=\mathbf{e}_i\mathbf{e}_j^\top+\mathbf{e}_j\mathbf{e}_i^\top.
\end{equation}
Thus~$M^{\{i,j\},f}$ is the martingale transform of~$M^f$ by~$A_{i,j}$.

On a jump of~$\mathcal{N}^k$ with sign~$\tau$, the increment of~$M^g$ is~$\tfrac12\widetilde{\nabla}^2_k\widetilde{g}+\tfrac\tau2\widetilde{\nabla}^0_k\widetilde{g}$. The increment of~$M^{\{i,j\},f}$ is the same expression after~$A_{i,j}$ has swapped the~$i$-th and~$j$-th coordinates, and is therefore supported only on jumps of~$\mathcal{N}^i$ and~$\mathcal{N}^j$. The product of these increments is a single one-sided pairing: using~\eqref{eq:grad02} and~$\tau^2=1$,
\begin{align*}
	\text{if the active coordinate is }k=i,
	 & \qquad
	\partial_j^{\tau}\widetilde{f}(T-t,\widetilde{X}_{t-})\,\partial_i^{\tau}\widetilde{g}(T-t,\widetilde{X}_{t-}), \\
	\text{if the active coordinate is }k=j,
	 & \qquad
	\partial_i^{\tau}\widetilde{f}(T-t,\widetilde{X}_{t-})\,\partial_j^{\tau}\widetilde{g}(T-t,\widetilde{X}_{t-}).
\end{align*}
Thus, by~\eqref{eq:qv-drivers},
\begin{align*}
	\mathrm{d}\bigl[M^{\{i,j\},f},M^g\bigr]_t
	 & =
	\partial_j^{\tau_{\mathcal{N}_t}}\widetilde{f}(T-t,\widetilde{X}_{t-})\,\partial_i^{\tau_{\mathcal{N}_t}}\widetilde{g}(T-t,\widetilde{X}_{t-})\,\mathrm{d}\mathcal{N}^i_t \\
	 & \qquad
	+
	\partial_i^{\tau_{\mathcal{N}_t}}\widetilde{f}(T-t,\widetilde{X}_{t-})\,\partial_j^{\tau_{\mathcal{N}_t}}\widetilde{g}(T-t,\widetilde{X}_{t-})\,\mathrm{d}\mathcal{N}^j_t.
\end{align*}
Since~$\tau_{\mathcal{N}_t}$ is independent of~$\mathcal{F}_{t-}$ and each~$\mathcal{N}^k$ has intensity~$2$, the freezing lemma yields
\begin{align*}
	\mathbb{E}_x\Bigl[\bigl[M^{\{i,j\},f},M^g\bigr]_T\Bigr]
	 & =\mathbb{E}_x\Biggl[
		                \int_0^T
		                \partial_j^{\tau_{\mathcal{N}_t}}\widetilde{f}(T-t,\widetilde{X}_{t-})\,\partial_i^{\tau_{\mathcal{N}_t}}\widetilde{g}(T-t,\widetilde{X}_{t-})\,\mathrm{d}\mathcal{N}^i_t\\
		                &\hspace{4em}
		                +\partial_i^{\tau_{\mathcal{N}_t}}\widetilde{f}(T-t,\widetilde{X}_{t-})\,\partial_j^{\tau_{\mathcal{N}_t}}\widetilde{g}(T-t,\widetilde{X}_{t-})\,\mathrm{d}\mathcal{N}^j_t
		                \Biggr]               \\
	 & =\mathbb{E}_x\,\mathbb{E}_x\Biggl[
		                              \int_0^T
		                              \partial_j^{\tau_{\mathcal{N}_t}}\widetilde{f}(T-t,\widetilde{X}_{t-})\,\partial_i^{\tau_{\mathcal{N}_t}}\widetilde{g}(T-t,\widetilde{X}_{t-})\,\mathrm{d}\mathcal{N}^i_t\\
		                              &\hspace{4em}
		                              +\partial_i^{\tau_{\mathcal{N}_t}}\widetilde{f}(T-t,\widetilde{X}_{t-})\,\partial_j^{\tau_{\mathcal{N}_t}}\widetilde{g}(T-t,\widetilde{X}_{t-})\,\mathrm{d}\mathcal{N}^j_t
		                              \Bigm|\mathcal{F}_{t-}
		                              \Biggr] \\
	 & =\mathbb{E}_x\int_0^T
	\Bigl(
	\partial_j^+\widetilde{f}\,\partial_i^+\widetilde{g}
	+\partial_j^-\widetilde{f}\,\partial_i^-\widetilde{g}
	+\partial_i^+\widetilde{f}\,\partial_j^+\widetilde{g}
	+\partial_i^-\widetilde{f}\,\partial_j^-\widetilde{g}
	\Bigr)(T-t,\widetilde{X}_{t-})\,\mathrm{d}t.
\end{align*}
Here~$\mathbb{P}_x$ denotes the law of~$(\widetilde{\mathcal{N}},\widetilde{X},(\tau_k))$ under the deterministic initial condition~$\widetilde{X}_0=x$, and~$\mathbb{E}_x$ is the corresponding expectation. The unqualified expectation
\begin{equation}
	\label{eq:unqualifiedE}
	\mathbb{E}
	:=\sum_{x\in\mathbb{Z}^N}\mathbb{E}_x
\end{equation}
includes the sum over all starting points. Averaging the hitting probabilities against the starting point gives~$\sum_x\mathbb{P}_x(\widetilde{X}_T=z)=1$ for every endpoint~$z$. Summing the equality above over the starting points~$x\in\mathbb{Z}^N$ of the random walk, the occupation identity
\(\sum_x\mathbb{E}_x H(\widetilde{X}_{t-})=\sum_y H(y)\) replaces the position of the walk by a free spatial variable~$y\). Writing that dummy as~$y\) to distinguish it from the starting point, we get
\begin{align*}
	\sum_{x\in\mathbb{Z}^N}
	\mathbb{E}_x\Bigl[\bigl[M^{\{i,j\},f},M^g\bigr]_T\Bigr]
	 & =\sum_{y\in\mathbb{Z}^N}\int_0^T
	\Bigl(
	\partial_j^+\widetilde{f}\,\partial_i^+\widetilde{g}
	+\partial_j^-\widetilde{f}\,\partial_i^-\widetilde{g}
	+\partial_i^+\widetilde{f}\,\partial_j^+\widetilde{g}
	+\partial_i^-\widetilde{f}\,\partial_j^-\widetilde{g}
	\Bigr)(T-t,y)\,\mathrm{d}t          \\
	 & =\sum_{y\in\mathbb{Z}^N}\int_0^T
	\Bigl(
	\bigl\langle A_{i,j}\,\widetilde{\nabla}^+\widetilde{f},\,\widetilde{\nabla}^+\widetilde{g}\bigr\rangle
	+\bigl\langle A_{i,j}\,\widetilde{\nabla}^-\widetilde{f},\,\widetilde{\nabla}^-\widetilde{g}\bigr\rangle
	\Bigr)(T-t,y)\,\mathrm{d}t.
\end{align*}
Here we have written
\[
	\widetilde{\nabla}^+\varphi
	=\bigl(\partial_1^+\varphi,\ldots,\partial_N^+\varphi\bigr),
	\qquad
	\widetilde{\nabla}^-\varphi
	=\bigl(\partial_1^-\varphi,\ldots,\partial_N^-\varphi\bigr)
\]
for the forward and backward discrete gradients, so that~$A_{i,j}$ acts as a quadratic form separately on each copy of~$\mathbb{R}^N$. The right-hand side is twice the truncated off-diagonal pairing~\eqref{eq:weak-off}. Letting~$T\to\infty$, one obtains
\[
	\langle R_{(i,j)}^2f,g\rangle
	=-\frac12\lim_{T\to\infty}
	\mathbb{E}\bigl[\bigl[M^{\{i,j\},f},M^g\bigr]_T\bigr].
\]

The same construction applies to a general coefficient matrix. Write~$A_{i,i}:=\mathbf{e}_i\mathbf{e}_i^\top$ for the diagonal elementary matrices, and assemble
\begin{equation}
	\label{eq:Aalpha}
	A_{\alpha,\mathcal{A}}
	=\sum_{i\in\mathcal{A}'}\alpha_i\,A_{i,i}
	+\sum_{\{i,j\}\in\mathcal{A}}\alpha_i\,A_{i,j}.
\end{equation}
The associated martingale transform starts at zero:
\begin{equation}
	\label{eq:Malpha-mart}
	M^{\alpha,\mathcal{A},f}_t
	=\frac12\int_0^t
	\Bigl(
	\bigl\langle A_{\alpha,\mathcal{A}}\widetilde{\nabla}^2\widetilde{f}(T-s,\widetilde{X}_{s-}),\,\mathrm{d}\bigl(\widetilde{\mathcal{N}}_s-2s\mathbf{1}\bigr)\bigr\rangle
	+\bigl\langle A_{\alpha,\mathcal{A}}\widetilde{\nabla}^0\widetilde{f}(T-s,\widetilde{X}_{s-}),\,\mathrm{d}\widetilde{X}_s\bigr\rangle
	\Bigr).
\end{equation}
Linearity of the quadratic covariation together with~\eqref{eq:RalphaA} yields
\[
	\langle R_{\alpha,\mathcal{A}}^2f,g\rangle
	=-\lim_{T\to\infty}
	\mathbb{E}\bigl[\bigl[M^{\alpha,\mathcal{A},f},M^g\bigr]_T\bigr].
\]
As in Gundy--Varopoulos, the pairing is realised by the conditional expectation of the transformed martingale given the endpoint of the background noise, using the unqualified expectation~\eqref{eq:unqualifiedE}. Under that mixture the averaged hitting probability of any site is~$1$, so there is no extra factor~$\mathbb{P}(\widetilde{X}_T=z)$.

\begin{lemma}
	\label{lem:GV}
	With the notation above, the second-order Riesz transform~\eqref{eq:RalphaA} is the Gundy--Varopoulos projection
	\[
		R_{\alpha,\mathcal{A}}^2f(x)
		=-\lim_{T\to\infty}
		\mathbb{E}\bigl[M^{\alpha,\mathcal{A},f}_T\bigm|\widetilde{X}_T=x\bigr].
	\]
	In particular~$B_{\mathfrak{i}}$ is the case~$N=2$,~$\mathcal{A}=\{(1,2),(2,1)\}$,~$\alpha=(-1,-1)$, and~$M^{\{1,2\},f}=-M^{\alpha,\mathcal{A},f}$ is the corresponding un-normalised transform.
\end{lemma}

\section{Dissipation of the Burkholder functional}
\label{sec:burkholder}

\subsection{The Burkholder functional}
\label{sec:U}

Write~$C_p:=p^\star-1$ and
\begin{equation}
	\label{eq:V}
	V(x,y)
	:=|y|^p-C_p^p\,|x|^p,
	\qquad
	x,y\in\mathbb{R}.
\end{equation}
The bound~$\|R_{\alpha,\mathcal{A}}^2\|_{\ell^p\to\ell^p}\le C_p$ is read off the martingale transform~$M^{\alpha,\mathcal{A},f}$ of Lemma~\ref{lem:GV}. Writing~$M^f$ for the heat-extension martingale of the input, the Gundy--Varopoulos projection and the \(\ell^p\) contractivity of the semigroup reduce the operator-norm estimate to
\[
	\mathbb{E}\,V\bigl(M^f_T,\,M^{\alpha,\mathcal{A},f}_T\bigr)\le 0
\]
in the large-height limit. The strategy is to dominate~$V$ by Burkholder's function~$U$ and to show that~$U$ is a supermartingale along the pair~$\bigl(M^f,M^{\alpha,\mathcal{A},f}\bigr)$, so that
\begin{equation}
	\label{eq:VU-chain}
	\mathbb{E}\,V\bigl(M^f_T,\,M^{\alpha,\mathcal{A},f}_T\bigr)
	\le
	\mathbb{E}\,U\bigl(M^f_T,\,M^{\alpha,\mathcal{A},f}_T\bigr)
	\le
	\mathbb{E}\,U\bigl(M^f_0,\,M^{\alpha,\mathcal{A},f}_0\bigr)
	\le 0.
\end{equation}
The first comparison is the pointwise inequality~$U\ge V$; the second is the dissipation of~$U$ at jumps, proved in the next subsections; the third holds because the transform starts at zero, so~$U(M^f_0,M^{\alpha,\mathcal{A},f}_0)=U(M^f_0,0)\le 0$ by~\eqref{eq:Uneg}. The function~$V$ itself does not have the concavity needed to run Itô's formula. Burkholder's substitute~\cite{Burkholder1984,Burkholder1991} is
\begin{equation}
	\label{eq:U}
	U(x,y)
	:=\alpha_p\bigl(|y|-C_p|x|\bigr)\bigl(|x|+|y|\bigr)^{p-1},
	\qquad
	\alpha_p
	:=p\Bigl(1-\frac1{p^\star}\Bigr)^{p-1}.
\end{equation}
It satisfies the three structural properties
\begin{align}
	\label{eq:UgeV}
	U(x,y) & \ge V(x,y)
	       &            & \text{for all }x,y,        \\
	\label{eq:Uneg}
	U(x,y) & \le 0
	       &            & \text{whenever }|y|\le|x|, \\
	\label{eq:U0}
	U(0,0) & =0.
\end{align}
In particular~$U(x,0)\le 0$. The function is not concave, but it is \emph{zigzag concave}: for every real~$\alpha$ with~$|\alpha|\le 1$, the restrictions
\begin{equation}
	\label{eq:zigzag}
	t\mapsto U(x+t,\,y+\alpha t),
	\qquad
	t\mapsto U(x+t,\,y-\alpha t)
\end{equation}
are concave. In other words~$U$ is concave along every line whose slope has absolute value at most~$1$. This is the analytic form of differential subordination, and it is the reason a direct Itô computation on a jump process is not immediate: the increments of~$\bigl(M^f,M^{\alpha,\mathcal{A},f}\bigr)$ are not constrained to a single line of slope at most~$1$ in absolute value. That point is taken up in Subsection~\ref{sec:remodel}.

Throughout the Itô computations we assume without loss of generality that~$U$ (and likewise~$V$) is of class~$C^2$, even though the expressions~\eqref{eq:U} and~\eqref{eq:V} fail to be smooth on the axes~$x=0$ and~$y=0$. This is arranged in the standard way of Burkholder, by adding a fictitious coordinate of size~$\varepsilon>0$ to the vectors~$(x,y)$. The moduli appearing in~$U$ and~$V$ are thereby replaced by smoothed Euclidean norms, all derivatives up to order two remain bounded on compact sets, and the structural properties~\eqref{eq:UgeV}--\eqref{eq:zigzag} persist uniformly as~$\varepsilon\to 0$.

\subsection{Dissipation of Burkholder's functional for the case~$N=2$}
\label{sec:dissipation}

We now specialise to dimension~$N=2$ and to the pair~$\{1,2\}$. Write~$A:=A_{1,2}$ for the elementary swap~\eqref{eq:Aij} and
\[
	M_t=M^f_t,
	\qquad
	M^A_t=M^{\{1,2\},f}_t
\]
for the heat-extension martingale of~$f$ and its transform by~$A$. (This is the un-normalised transform attached to~$B_{\mathfrak{i}}=-2R_{(1,2)}^2\).) Thus~$M_t=\widetilde{f}(T-t,\widetilde{X}_t)$ with initial value~$M_0=\widetilde{f}(T,\widetilde{X}_0)$, while the transform starts at zero:
	\[
	M^A_0=0.
\]
Using the Itô formula for càdlàg processes we obtain
\begin{equation}
	\label{eq:itoU}
	\begin{aligned}
		\mathbb{E}\,U(M_t,M^A_t)
		 & =\mathbb{E}\,U(M_0,M^A_0) \\
		 & \qquad
		+\mathbb{E}\sum_{0<s\le t}
		\Bigl(
		U(M_s,M^A_s)-U(M_{s-},M^A_{s-})
		-\nabla U(M_{s-},M^A_{s-})\cdot\Delta(M_s,M^A_s)
		\Bigr)                       \\
		 & =\mathbb{E}\,U(M_0,M^A_0) \\
		 & \qquad
		+\mathbb{E}\,\mathbb{E}\Biggl(
		\sum_{0<s\le t}
		\Bigl(
		U(M_s,M^A_s)-U(M_{s-},M^A_{s-})
		-\nabla U(M_{s-},M^A_{s-})\cdot\Delta(M_s,M^A_s)
		\Bigr)
		\Bigm|\mathcal{F}_{s-}
		\Biggr).
	\end{aligned}
\end{equation}
The processes~$M$ and~$M^A$ are martingales, so the first-order continuous compensator of the jumps cancels, in expectation, the first-order contribution of the jumps themselves; only the second-order jump remainder survives. The inner conditional increment is the dissipation at the jump time~$s$. We write
\begin{equation}
	\label{eq:Ds}
	\mathcal{D}_s
	:=\mathbb{E}\Biggl(
	U(M_s,M^A_s)-U(M_{s-},M^A_{s-})
	-\nabla U(M_{s-},M^A_{s-})\cdot\Delta(M_s,M^A_s)
	\Bigm|\mathcal{F}_{s-}
	\Biggr)
\end{equation}
and
\[
	\mathbb{E}\,U(M_t,M^A_t)
	=\mathbb{E}\,U(M_0,M^A_0)
	+\mathbb{E}\sum_{0<s\le t}\mathcal{D}_s.
\]
Following Burkholder, the whole estimate reduces to the sign of this dissipation.

\begin{lemma}
	\label{lem:dissipation}
	For every jump time~$s\in(0,T]$ one has~$\mathcal{D}_s\le 0$.
\end{lemma}

The proof is given in Subsection~\ref{sec:remodel}. We first record the consequence for~$B_{\mathfrak{i}}$.

\subsection{The sharp~$\ell^p$ bound for~$B_{\mathfrak{i}}$}
\label{sec:sharpBi}

Assume Lemma~\ref{lem:dissipation}. Then~\eqref{eq:itoU} yields the supermartingale inequality
\begin{equation}
	\label{eq:super}
	\mathbb{E}\,U(M_t,M^A_t)
	\le
	\mathbb{E}\,U(M_0,M^A_0),
	\qquad
	0\le t\le T.
\end{equation}
Since~$M^A_0=0$ and~$U(x,0)\le 0$ by~\eqref{eq:Uneg}, the right-hand side is non-positive. Combined with~$U\ge V$ this gives
\[
	\mathbb{E}\,V(M_T,M^A_T)
	\le
	\mathbb{E}\,U(M_T,M^A_T)
	\le 0,
\]
that is
\begin{equation}
	\label{eq:moment}
	\mathbb{E}\,|M^A_T|^p
	\le
	C_p^p\,\mathbb{E}\,|M_T|^p.
\end{equation}
It remains to pass from these martingale moments to the~$\ell^p(\mathbb{Z}^2)$ norms of~$f$ and~$B_{\mathfrak{i}}f$, following the averaging procedure of Arcozzi--Domelevo--Petermichl~\cite{ArcDomPet2016a}.

Start the background noise independently at every lattice point, and write~$\mathbb{E}_x$ for the law with~$\widetilde{X}_0=x$. Summing~\eqref{eq:moment} over the starting points produces the deterministic inequality
\begin{equation}
	\label{eq:summed}
	\sum_{x\in\mathbb{Z}^2}
	\mathbb{E}_x\,|M^A_T|^p
	\le
	C_p^p
	\sum_{x\in\mathbb{Z}^2}
	\mathbb{E}_x\,|M_T|^p.
\end{equation}
On the other hand, Lemma~\ref{lem:GV} and Jensen's inequality for the conditional expectation give
\begin{align*}
	\|B_{\mathfrak{i}}f\|_{\ell^p(\mathbb{Z}^2)}^p
	 & =\sum_{z\in\mathbb{Z}^2}
	\Bigl|
	\lim_{T\to\infty}
	\mathbb{E}\bigl[M^A_T\bigm|\widetilde{X}_T=z\bigr]
	\Bigr|^p
	\le
	\liminf_{T\to\infty}
	\sum_{z\in\mathbb{Z}^2}
	\mathbb{E}\bigl[|M^A_T|^p\bigm|\widetilde{X}_T=z\bigr] \\
	 & =\liminf_{T\to\infty}
	\mathbb{E}\,|M^A_T|^p
	=\liminf_{T\to\infty}
	\sum_{x\in\mathbb{Z}^2}
	\mathbb{E}_x\,|M^A_T|^p,
\end{align*}
where the last step uses that the law of the endpoint, averaged over starting points, is translation invariant on~$\mathbb{Z}^2$. For the right-hand side of~\eqref{eq:summed} one has~$M_T=M^f_T=\widetilde{f}(0,\widetilde{X}_T)=f(\widetilde{X}_T)$. Averaging over starting points as in~\eqref{eq:unqualifiedE} therefore yields
\[
	\sum_{x\in\mathbb{Z}^2}
	\mathbb{E}_x\,|M_T|^p
	=\|f\|_{\ell^p(\mathbb{Z}^2)}^p.
\] Combining these comparisons and sending~$T\to\infty$ in~\eqref{eq:summed} gives
\[
	\|B_{\mathfrak{i}}f\|_{\ell^p(\mathbb{Z}^2)}
	\le
	C_p\,\|f\|_{\ell^p(\mathbb{Z}^2)}
	= (p^\star-1)\,\|f\|_{\ell^p(\mathbb{Z}^2)}.
\]
This is the upper bound of Theorem~\ref{thm:imag}. The comparison~$\mathbb{E}\,U(M_T,M^A_T)\le\mathbb{E}\,U(M_0,M^A_0)$ used above rests on the dissipation~$\mathcal{D}_s\le 0$, whose second summand is non-positive by the zigzag concavity~\eqref{eq:zigzag}. The matching lower bound is inherited from the continuous operator of Geiss--Montgomery-Smith--Saksman~\cite{GMSS2009} by sending the mesh size to zero.

\subsection{A remodeling trick and the proof of Lemma~\ref{lem:dissipation}}
\label{sec:remodel}

In order to estimate the dissipation of the Burkholder functional we observe that the dissipation does not only involve the graph Laplacian of~$U$ but also a drift term. Added to that, the points involved in this graph Laplacian do not produce martingale jumps, and can therefore not be estimated directly by Burkholder's original argument. Finally, as mentioned above, the jumps of the martingale and of its transform are not differentially subordinate.

Two ingredients are needed. The first concerns the drift. The dissipation is of the form
\[
	\mathcal{D}
	=\frac14\sum_{i=1}^4
	\Bigl(U(v_i)-U(v_0)-U'(v_0)\cdot(v_i-v_0)\Bigr),
\]
with non-martingale jumps averaging to
\[
	a
	:=\frac14\sum_{i=1}^4(v_i-v_0).
\]
That is an expression
\[
	\frac14\sum_{i=1}^4 U(v_i)-U(v_0)-U'(v_0)\cdot a,
\]
and the first-order approximation~$U(v_0)+U'(v_0)\cdot a\simeq U(v_0+a)$ suggests that the dissipation should not be far from
\[
	\frac14\sum_{i=1}^4 U(v_i)-U(v_0+a).
\]
Making explicit that the jumps should now be read from~$v_0+a$ to~$v_i$, and writing the error in the approximation above, the first step is to split
\begin{equation}
	\label{eq:Dsplit-model}
	\begin{aligned}
		\mathcal{D}
		 & =\frac14\sum_{i=1}^4
		\Bigl(U\bigl((v_0+a)+(v_i-v_0-a)\bigr)-U(v_0+a)\Bigr)
		+\Bigl(U(v_0+a)-U(v_0)-U'(v_0)\cdot a\Bigr) \\
		 & =:\mathcal{D}_1+\mathcal{D}_2.
	\end{aligned}
\end{equation}
This identity has to be adapted to the present transform. The crucial structural fact is that the \emph{same} shift~$a$ appears for~$M$ and for~$M^A$. That common drift lets us estimate~$\mathcal{D}_2$ from the zigzag concavity of~$U$ along the diagonal.

The second ingredient concerns~$\mathcal{D}_1$ and the fact that one has only a weak form of differential subordination of~$M^A$ with respect to~$M$. To handle this we remodel the two martingales by inserting intermediate points. This is again possible because of the specific form of the transform~$A$. The details follow.

Fix a jump time~$s$. Write~$M$ and~$M^A$ for the left limits~$M_{s-}$ and~$M^A_{s-}$, and write~$\partial_1U$,~$\partial_2U$ for the partial derivatives of~$U$ in its two arguments. Conditioning on~$\mathcal{F}_{s-}$ averages the four possible jumps of~$(\mathcal{N}^1,\mathcal{N}^2)$ with signs~$\tau=\pm$, each of intensity~$1$, and yields
\begin{equation}
	\label{eq:Ds-avg}
	\mathcal{D}_s
	=\frac14\sum_{k=1}^2\sum_{\tau=\pm}
	\Bigl(
	U\bigl(M+\delta^s_{k,\tau},\,M^A+\delta^{A,s}_{k,\tau}\bigr)
	-U(M,M^A)
	-\partial_1U(M,M^A)\,\delta^s_{k,\tau}
	-\partial_2U(M,M^A)\,\delta^{A,s}_{k,\tau}
	\Bigr).
\end{equation}
The increments of~$M$ are
\begin{align*}
	\delta^s_{1,+}
	 & =\frac12\bigl(\partial_1^2\widetilde{f}+\partial_1^0\widetilde{f}\bigr)(T-s,\widetilde{X}_{s-})
	=\partial_1^+\widetilde{f}(T-s,\widetilde{X}_{s-}),                                                \\
	\delta^s_{1,-}
	 & =\frac12\bigl(\partial_1^2\widetilde{f}-\partial_1^0\widetilde{f}\bigr)(T-s,\widetilde{X}_{s-})
	=-\partial_1^-\widetilde{f}(T-s,\widetilde{X}_{s-}),                                               \\
	\delta^s_{2,+}
	 & =\frac12\bigl(\partial_2^2\widetilde{f}+\partial_2^0\widetilde{f}\bigr)(T-s,\widetilde{X}_{s-})
	=\partial_2^+\widetilde{f}(T-s,\widetilde{X}_{s-}),                                                \\
	\delta^s_{2,-}
	 & =\frac12\bigl(\partial_2^2\widetilde{f}-\partial_2^0\widetilde{f}\bigr)(T-s,\widetilde{X}_{s-})
	=-\partial_2^-\widetilde{f}(T-s,\widetilde{X}_{s-}).
\end{align*}
The increments of~$M^A$ are those of~$M$ with the two coordinates swapped,
\begin{align*}
	\delta^{A,s}_{1,+} & =\delta^s_{2,+}, &
	\delta^{A,s}_{1,-} & =\delta^s_{2,-}, &
	\delta^{A,s}_{2,+} & =\delta^s_{1,+}, &
	\delta^{A,s}_{2,-} & =\delta^s_{1,-}.
\end{align*}
The mean increment is therefore the same for both martingales:
\begin{equation}
	\label{eq:a-same}
	a
	:=\frac14\sum_{k=1}^2\sum_{\tau=\pm}\delta^s_{k,\tau}
	=\frac14\sum_{k=1}^2\sum_{\tau=\pm}\delta^{A,s}_{k,\tau},
\end{equation}
and consequently
\[
	\frac14\sum_{k=1}^2\sum_{\tau=\pm}(\delta^s_{k,\tau}-a)
	=\frac14\sum_{k=1}^2\sum_{\tau=\pm}(\delta^{A,s}_{k,\tau}-a)
	=0.
\]
Unwinding the one-sided differences,
\[
	\partial_k^+\widetilde{f}-\partial_k^-\widetilde{f}
	=\widetilde{f}(\,\cdot\,+e_k)+\widetilde{f}(\,\cdot\,-e_k)-2\widetilde{f}
	=\partial_k^2\widetilde{f},
\]
the sum of the four increments is the spatial Laplacian, and therefore
\begin{equation}
	\label{eq:a-lap}
	a
	=\frac14\sum_{k=1}^2\sum_{\tau=\pm}\delta^s_{k,\tau}
	=\frac14\Delta\widetilde{f}(T-s,\widetilde{X}_{s-}).
\end{equation}
This common average increment is the drift term advertised above. Substituting~\eqref{eq:a-same} into~\eqref{eq:Ds-avg} and rearranging as in~\eqref{eq:Dsplit-model} gives
\begin{equation}
	\label{eq:D12}
	\begin{aligned}
		\mathcal{D}_s
		 & =\frac14\sum_{k=1}^2\sum_{\tau=\pm}
		\Bigl(
		U\bigl((M+a)+(\delta^s_{k,\tau}-a),\,(M^A+a)+(\delta^{A,s}_{k,\tau}-a)\bigr)
		-U(M+a,M^A+a)
		\Bigr)                                    \\
		 & \qquad
		+U(M+a,M^A+a)-U(M,M^A)
		-\partial_1U(M,M^A)\,a
		-\partial_2U(M,M^A)\,a                    \\
		 & =:\mathcal{D}_{s,1}+\mathcal{D}_{s,2}.
	\end{aligned}
\end{equation}
The second summand is the second-order remainder of~$U$ along the diagonal increment~$(a,a)$. Zigzag concavity~\eqref{eq:zigzag} of~$U$ along lines of slope~$+1$ therefore yields
\begin{equation}
	\label{eq:D2neg}
	\mathcal{D}_{s,2}\le 0.
\end{equation}

It remains to show~$\mathcal{D}_{s,1}\le 0$. Introduce the remodeled increments
\begin{equation}
	\label{eq:eta}
	\eta
	:=\frac12\bigl((\delta^s_{2,-}-a)+(\delta^s_{1,-}-a)\bigr),
	\qquad
	\eta^A
	:=\frac12\bigl((\delta^{A,s}_{2,-}-a)+(\delta^{A,s}_{1,-}-a)\bigr).
\end{equation}
The identity~$\sum_{k,\tau}(\delta^s_{k,\tau}-a)=0$ gives~$-\eta=\tfrac12\bigl((\delta^s_{2,+}-a)+(\delta^s_{1,+}-a)\bigr)$, and likewise for the transformed increments. The swap relations~$\delta^{A,s}_{1,\pm}=\delta^s_{2,\pm}$ and~$\delta^{A,s}_{2,\pm}=\delta^s_{1,\pm}$ force
\[
	\eta^A=\eta,
	\qquad
	-\eta^A=-\eta.
\]
Thus the two minus-jumps are centered at the diagonal point~$(M+a+\eta,\,M^A+a+\eta)$, and the two plus-jumps are centered at~$(M+a-\eta,\,M^A+a-\eta)$. Inserting these intermediate values,
\begin{align*}
	\mathcal{D}_{s,1}
	 & =\frac14\sum_{k=1}^2\sum_{\tau=\pm}
	U\bigl((M+a)+(\delta^s_{k,\tau}-a),\,(M^A+a)+(\delta^{A,s}_{k,\tau}-a)\bigr)
	-U(M+a,M^A+a)                          \\
	 & =\Biggl(
	\frac14\sum_{k,\tau}U(\cdots)
	-\frac12 U(M+a+\eta,\,M^A+a+\eta)
	-\frac12 U(M+a-\eta,\,M^A+a-\eta)
	\Biggr)                                \\
	 & \qquad
	+\Biggl(
	\frac12 U(M+a+\eta,\,M^A+a+\eta)
	+\frac12 U(M+a-\eta,\,M^A+a-\eta)
	-U(M+a,M^A+a)
	\Biggr).
\end{align*}
The second parenthesis is again a diagonal second-order remainder, hence non-positive by zigzag concavity along slope~$+1$. The first parenthesis splits according to the sign of~$\tau$:
\begin{align*}
	 & \frac14\sum_{k=1}^2
	U\bigl((M+a)+(\delta^s_{k,-}-a),\,(M^A+a)+(\delta^{A,s}_{k,-}-a)\bigr)
	-\frac12 U(M+a+\eta,\,M^A+a+\eta) \\
	 & \qquad
	+\frac14\sum_{k=1}^2
	U\bigl((M+a)+(\delta^s_{k,+}-a),\,(M^A+a)+(\delta^{A,s}_{k,+}-a)\bigr)
	-\frac12 U(M+a-\eta,\,M^A+a-\eta).
\end{align*}
For the minus pair, the two deviations from the midpoint~$(M+a+\eta,\,M^A+a+\eta)$ are opposite and of the form~$(h,-h)$. Indeed,
\[
	(\delta^s_{1,-}-a)-\eta
	=\frac12\bigl((\delta^s_{1,-}-a)-(\delta^s_{2,-}-a)\bigr)
	=-
	\bigl((\delta^{A,s}_{1,-}-a)-\eta\bigr),
\]
so the chord joining the two minus-images has slope~$-1$. Zigzag concavity along lines of slope~$-1$ therefore makes the minus contribution non-positive. The plus pair is symmetric about~$(M+a-\eta,\,M^A+a-\eta)$ with the same slope~$-1$, and is likewise non-positive. Hence~$\mathcal{D}_{s,1}\le 0$.

Combining with~\eqref{eq:D2neg} we obtain~$\mathcal{D}_s\le 0$, which is Lemma~\ref{lem:dissipation}. The estimate on~$B_{\mathfrak{i}}$ then follows as in Subsection~\ref{sec:sharpBi}.

\subsection{The sharp~$\ell^p$ bound for~$R_{\alpha,\mathcal{A}}^2$}
\label{sec:sharp-general}

The same chain applies to a general operator~\eqref{eq:RalphaA}. Let~$M^f$ be the heat-extension martingale of the input and let~$M^{\alpha,\mathcal{A},f}$ be the transform~\eqref{eq:Malpha-mart} by the assembled matrix~$A_{\alpha,\mathcal{A}}$ of~\eqref{eq:Aalpha}. Lemma~\ref{lem:GV} realises~$R_{\alpha,\mathcal{A}}^2f$ as the Gundy--Varopoulos projection of~$-M^{\alpha,\mathcal{A},f}_T$. The comparison~\eqref{eq:VU-chain} therefore reduces the bound
\[
	\bigl\|R_{\alpha,\mathcal{A}}^2\bigr\|_{\ell^p(\mathbb{Z}^N)\to\ell^p(\mathbb{Z}^N)}
	\le p^\star-1
\]
to the dissipation of~$U$ along the pair~$\bigl(M^f,M^{\alpha,\mathcal{A},f}\bigr)$.

Each elementary block of~$A_{\alpha,\mathcal{A}}$ is either a coordinate projection~$A_{i,i}$ or a transposition~$A_{i,j}$ of two coordinates.
We want to use this structure to control again the dissipation of Burkholder's functional. For that, we remodel the~$2N$ jumps according to the elementary block structure of~$A_{\alpha,\mathcal{A}}$. Let us first rewrite the Itô formula after~\eqref{eq:Ds}, now for the pair~$\bigl(M^f,M^{\alpha,\mathcal{A},f}\bigr)$:
\begin{equation}
	\label{eq:itoU-general}
	\begin{aligned}
		\mathbb{E}\,U\bigl(M^f_t,M^{\alpha,\mathcal{A},f}_t\bigr)
		 & =\mathbb{E}\,U\bigl(M^f_0,M^{\alpha,\mathcal{A},f}_0\bigr) \\
		 & \qquad
		+\mathbb{E}\sum_{0<s\le t}
		\mathcal{D}_s,
	\end{aligned}
\end{equation}
where the dissipation~$\mathcal{D}_s$ now involves all directions of the augmented gradients and splits into diagonal and off-diagonal blocks,
\begin{equation}
	\label{eq:Ds-general}
	\mathcal{D}_s
	=\sum_{i\in\mathcal{A}'}\mathcal{D}_s^{(i)}
	+\sum_{\{i,j\}\in\mathcal{A}}\mathcal{D}_s^{(i,j)}.
\end{equation}
Write
\[
	R_s
	:=
	U\bigl(M^f_s,M^{\alpha,\mathcal{A},f}_s\bigr)
	-U\bigl(M^f_{s-},M^{\alpha,\mathcal{A},f}_{s-}\bigr)
	-\nabla U\bigl(M^f_{s-},M^{\alpha,\mathcal{A},f}_{s-}\bigr)
	\cdot
	\bigl(\Delta M^f_s,\,\Delta M^{\alpha,\mathcal{A},f}_s\bigr)
\]
for the Itô remainder at time~$s$. For an unmatched index~$i\in\mathcal{A}'$ only the two jumps of~$\mathcal{N}^i$ fire, with transformed increment scaled by~$\alpha_i$:
\begin{equation}
	\label{eq:Ds-diag}
	\mathcal{D}_s^{(i)}
	=\mathbb{E}\Biggl(
	\sum_{\tau=\pm}
	R_s\,
	\mathbf{1}_{\{\Delta\mathcal{N}^i_s=1,\,\tau_{\mathcal{N}_s}=\tau\}}
	\Bigm|\mathcal{F}_{s-}
	\Biggr).
\end{equation}
For a pair~$\{i,j\}\in\mathcal{A}$ the four jumps of~$\mathcal{N}^i$ and~$\mathcal{N}^j$ fire, and~$A_{i,j}$ swaps the two coordinates before the factor~$\alpha_i$:
\begin{equation}
	\label{eq:Ds-off}
	\mathcal{D}_s^{(i,j)}
	=\mathbb{E}\Biggl(
	\sum_{k\in\{i,j\}}
	\sum_{\tau=\pm}
	R_s\,
	\mathbf{1}_{\{\Delta\mathcal{N}^k_s=1,\,\tau_{\mathcal{N}_s}=\tau\}}
	\Bigm|\mathcal{F}_{s-}
	\Biggr).
\end{equation}

Fix a jump time~$s$ and write~$M:=M^f_{s-}$,~$M^\alpha:=M^{\alpha,\mathcal{A},f}_{s-}$ for the left limits. The possible increments of~$M^f$ are, as between~\eqref{eq:Ds-avg} and~\eqref{eq:a-same},
\begin{align*}
	\delta^s_{k,+}
	 & =\partial_k^+\widetilde{f}(T-s,\widetilde{X}_{s-}),
	\qquad
	\delta^s_{k,-}
	=-\partial_k^-\widetilde{f}(T-s,\widetilde{X}_{s-}),
	\qquad
	k=1,\ldots,N.
\end{align*}
The increments of the transform are obtained by applying~$A_{\alpha,\mathcal{A}}$: on a jump of coordinate~$k$ with sign~$\tau$,
\[
	\delta^{\alpha,s}_{k,\tau}
	=\bigl(A_{\alpha,\mathcal{A}}\delta^s_{\,\cdot\,,\tau}\bigr)_k.
\]
In particular, for an unmatched index~$i\in\mathcal{A}'$ one has~$\delta^{\alpha,s}_{i,\tau}=\alpha_i\delta^s_{i,\tau}$, while for a pair~$\{i,j\}\in\mathcal{A}$ the transposition~$A_{i,j}$ yields
\[
	\delta^{\alpha,s}_{i,\tau}=\alpha_i\delta^s_{j,\tau},
	\qquad
	\delta^{\alpha,s}_{j,\tau}=\alpha_i\delta^s_{i,\tau}.
\]
Each of the~$2N$ signed jumps has intensity~$1$. Conditioning on~$\mathcal{F}_{s-}$ therefore rewrites~\eqref{eq:Ds-diag} and~\eqref{eq:Ds-off} as
\begin{equation}
	\label{eq:Ds-diag-avg}
	\begin{aligned}
		\mathcal{D}_s^{(i)}
		 & =\frac12\sum_{\tau=\pm}
		\Bigl(
		U\bigl(M+\delta^s_{i,\tau},\,M^\alpha+\delta^{\alpha,s}_{i,\tau}\bigr)
		-U(M,M^\alpha)
		-\partial_1U(M,M^\alpha)\,\delta^s_{i,\tau}
		-\partial_2U(M,M^\alpha)\,\delta^{\alpha,s}_{i,\tau}
		\Bigr)
	\end{aligned}
\end{equation}
and
\begin{equation}
	\label{eq:Ds-off-avg}
	\begin{aligned}
		\mathcal{D}_s^{(i,j)}
		 & =\frac14\sum_{k\in\{i,j\}}\sum_{\tau=\pm}
		\Bigl(
		U\bigl(M+\delta^s_{k,\tau},\,M^\alpha+\delta^{\alpha,s}_{k,\tau}\bigr)
		-U(M,M^\alpha)                               \\
		 & \hspace{4.2em}
		-\partial_1U(M,M^\alpha)\,\delta^s_{k,\tau}
		-\partial_2U(M,M^\alpha)\,\delta^{\alpha,s}_{k,\tau}
		\Bigr).
	\end{aligned}
\end{equation}
The mean increment of each block is the same for both coordinates of the pair~$(M,M^\alpha)$ up to the coefficient~$\alpha$. For an unmatched index we set
\begin{equation}
	\label{eq:ai}
	a_i
	:=\frac12\sum_{\tau=\pm}\delta^s_{i,\tau}
	=\frac12\partial_i^2\widetilde{f}(T-s,\widetilde{X}_{s-}),
\end{equation}
so that~$\tfrac12\sum_{\tau}(\delta^s_{i,\tau}-a_i)=0$ and~$\tfrac12\sum_{\tau}(\delta^{\alpha,s}_{i,\tau}-\alpha_i a_i)=0$. For a pair we set
\begin{equation}
	\label{eq:aij}
	a_{i,j}
	:=\frac14\sum_{k\in\{i,j\}}\sum_{\tau=\pm}\delta^s_{k,\tau}
	=\frac14\bigl(\partial_i^2+\partial_j^2\bigr)\widetilde{f}(T-s,\widetilde{X}_{s-}),
\end{equation}
so that the four centered increments of~$M$ and of~$M^\alpha$ both average to zero, the latter after the factor~$\alpha_i$.

Splitting each block as in~\eqref{eq:D12} therefore gives
\begin{equation}
	\label{eq:D12-diag}
	\begin{aligned}
		\mathcal{D}_s^{(i)}
		 & =\frac12\sum_{\tau=\pm}
		\Bigl(
		U\bigl(M+\delta^s_{i,\tau},\,M^\alpha+\delta^{\alpha,s}_{i,\tau}\bigr)
		-U(M+a_i,M^\alpha+\alpha_i a_i)
		\Bigr)                                                \\
		 & \qquad
		+U(M+a_i,M^\alpha+\alpha_i a_i)-U(M,M^\alpha)
		-\partial_1U(M,M^\alpha)\,a_i
		-\partial_2U(M,M^\alpha)\,\alpha_i a_i                \\
		 & =:\mathcal{D}_{s,1}^{(i)}+\mathcal{D}_{s,2}^{(i)}.
	\end{aligned}
\end{equation}
and
\begin{equation}
	\label{eq:D12-off}
	\begin{aligned}
		\mathcal{D}_s^{(i,j)}
		 & =\frac14\sum_{k\in\{i,j\}}\sum_{\tau=\pm}
		\Bigl(
		U\bigl(M+\delta^s_{k,\tau},\,M^\alpha+\delta^{\alpha,s}_{k,\tau}\bigr)
		-U(M+a_{i,j},M^\alpha+\alpha_i a_{i,j})
		\Bigr)                                                    \\
		 & \qquad
		+U(M+a_{i,j},M^\alpha+\alpha_i a_{i,j})-U(M,M^\alpha)
		-\partial_1U(M,M^\alpha)\,a_{i,j}
		-\partial_2U(M,M^\alpha)\,\alpha_i a_{i,j}                \\
		 & =:\mathcal{D}_{s,1}^{(i,j)}+\mathcal{D}_{s,2}^{(i,j)}.
	\end{aligned}
\end{equation}
The terms~$\mathcal{D}_{s,2}^{(i)}$ and~$\mathcal{D}_{s,2}^{(i,j)}$ are the second-order remainders of~$U$ along the block shifts~$(a_i,\alpha_i a_i)$ and~$(a_{i,j},\alpha_i a_{i,j})$. Both are therefore non-positive by the zigzag concavity~\eqref{eq:zigzag}: each shift has slope~$\alpha_i$ with~$|\alpha_i|\le 1$, which is the same property used before~\eqref{eq:D2neg} and at the end of Subsection~\ref{sec:sharpBi}.

Likewise,~$\mathcal{D}_{s,1}^{(i)}$ is a two-point average. The two remodeled increments from the block mean~$(M+a_i,\,M^\alpha+\alpha_i a_i)$ are opposite and of the form~$(h,\alpha_i h)$. The chord therefore has slope~$\alpha_i$, and two-point concavity along a line of slope at most~$1$ in absolute value---again~\eqref{eq:zigzag}---yields~$\mathcal{D}_{s,1}^{(i)}\le 0$.

It remains to treat~$\mathcal{D}_{s,1}^{(i,j)}$. Introduce intermediate points in the pair~$\{i,j\}$ exactly as in~\eqref{eq:eta}: write
\begin{equation}
	\label{eq:eta-ij}
	\eta_{i,j}
	:=\frac12\bigl((\delta^s_{j,-}-a_{i,j})+(\delta^s_{i,-}-a_{i,j})\bigr)
\end{equation}
for the midpoint of the two minus-increments.

We can in fact argue that the dissipation has the correct sign without calculation, by comparing with the proof of Subsection~\ref{sec:remodel}. Recall that the remodeling trick consists in introducing fictitious intermediate points and corresponding intermediate dyadic martingale jumps, observing that all of those jumps have a small enough slope. In the present case, if one first pretends that~$\alpha_i=\pm 1$, the configuration reduces to the four-point picture of Subsection~\ref{sec:remodel}: the two minus-jumps are centered at a fictitious point on the diagonal of slope~$+1$, the two plus-jumps at the opposite fictitious point, and every chord through those intermediate points has slope~$\pm 1$. In particular all jumps, including those that only exist after remodeling, have slope of absolute value~$1$. The actual coefficient~$\alpha_i$ is obtained from this picture by a linear contraction (or reflection, if~$\alpha_i<0$) of factor~$\alpha_i$ in the~$y$-direction, centered at the left-limit point~$(M,M^\alpha)$. This map sends a chord of slope~$\pm 1$ to a chord of slope~$\pm\alpha_i$. All four remodeled jumps, and the two intermediate points, therefore have slope of absolute value~$|\alpha_i|\le 1$. Zigzag concavity~\eqref{eq:zigzag} applies on each of those lines, and the same splitting as after~\eqref{eq:eta} yields~$\mathcal{D}_{s,1}^{(i,j)}\le 0$.

Let us now put this remark into equations. We follow the steps of~\eqref{eq:eta} and further. Write~$a:=a_{i,j}$ for the block mean~\eqref{eq:aij}. The identity~$\sum_{k\in\{i,j\}}\sum_{\tau}(\delta^s_{k,\tau}-a)=0$ gives
\[
	-\eta_{i,j}
	=\frac12\bigl((\delta^s_{j,+}-a)+(\delta^s_{i,+}-a)\bigr),
\]
and the transform increments satisfy~$\delta^{\alpha,s}_{i,\tau}=\alpha_i\delta^s_{j,\tau}$,~$\delta^{\alpha,s}_{j,\tau}=\alpha_i\delta^s_{i,\tau}$, hence
\begin{equation}
	\label{eq:eta-alpha}
	\eta^\alpha_{i,j}
	:=\frac12\bigl((\delta^{\alpha,s}_{j,-}-\alpha_i a)+(\delta^{\alpha,s}_{i,-}-\alpha_i a)\bigr)
	=\alpha_i\eta_{i,j}.
\end{equation}
Thus the two minus-jumps are centered at
\[
	\bigl(M+a+\eta_{i,j},\,M^\alpha+\alpha_i a+\alpha_i\eta_{i,j}\bigr),
\]
and the two plus-jumps at
\[
	\bigl(M+a-\eta_{i,j},\,M^\alpha+\alpha_i a-\alpha_i\eta_{i,j}\bigr).
\]
Inserting these intermediate values as after~\eqref{eq:eta},
\begin{align*}
	\mathcal{D}_{s,1}^{(i,j)}
	 & =\frac14\sum_{k\in\{i,j\}}\sum_{\tau=\pm}
	U\bigl((M+a)+(\delta^s_{k,\tau}-a),\,(M^\alpha+\alpha_i a)+(\delta^{\alpha,s}_{k,\tau}-\alpha_i a)\bigr) \\
	 & \qquad
	-U(M+a,M^\alpha+\alpha_i a)                                                                              \\
	 & =\Biggl(
	\frac14\sum_{k,\tau}U(\cdots)
	-\frac12 U(M+a+\eta_{i,j},\,M^\alpha+\alpha_i a+\alpha_i\eta_{i,j})                                      \\
	 & \qquad
	-\frac12 U(M+a-\eta_{i,j},\,M^\alpha+\alpha_i a-\alpha_i\eta_{i,j})
	\Biggr)                                                                                                  \\
	 & \qquad
	+\Biggl(
	\frac12 U(M+a+\eta_{i,j},\,M^\alpha+\alpha_i a+\alpha_i\eta_{i,j})                                       \\
	 & \qquad
	+\frac12 U(M+a-\eta_{i,j},\,M^\alpha+\alpha_i a-\alpha_i\eta_{i,j})
	-U(M+a,M^\alpha+\alpha_i a)
	\Biggr).
\end{align*}
The second parenthesis is the second-order remainder of~$U$ along the increment~$(\eta_{i,j},\alpha_i\eta_{i,j})$, hence non-positive by~\eqref{eq:zigzag}. The first parenthesis splits according to the sign of~$\tau$. For the minus pair the two deviations from the midpoint
\(\bigl(M+a+\eta_{i,j},\,M^\alpha+\alpha_i a+\alpha_i\eta_{i,j}\bigr)\)
are opposite and of the form~$(h,-\alpha_i h)$. Indeed,
\begin{align*}
	(\delta^s_{i,-}-a)-\eta_{i,j}
	 & =\frac12\bigl((\delta^s_{i,-}-a)-(\delta^s_{j,-}-a)\bigr), \\
	(\delta^{\alpha,s}_{i,-}-\alpha_i a)-\alpha_i\eta_{i,j}
	 & =\alpha_i\bigl((\delta^s_{j,-}-a)-\eta_{i,j}\bigr)
	=-\alpha_i\bigl((\delta^s_{i,-}-a)-\eta_{i,j}\bigr),
\end{align*}
so the chord joining the two minus-images has slope~$-\alpha_i$. Zigzag concavity~\eqref{eq:zigzag} along that line makes the minus contribution non-positive. The plus pair is symmetric about
\(\bigl(M+a-\eta_{i,j},\,M^\alpha+\alpha_i a-\alpha_i\eta_{i,j}\bigr)\)
with the same slope~$-\alpha_i$, and is likewise non-positive. Hence~$\mathcal{D}_{s,1}^{(i,j)}\le 0$.

Combining the four signs we obtain~$\mathcal{D}_s\le 0$. The comparison~\eqref{eq:VU-chain} follows, and Theorem~\ref{thm:main} is proved.

\section*{AI disclosure statement}
We have used Grok 4.6, built by xAI, to accelerate the write-up, bibliography searches and to check for correctness or
typos. All the ideas, models and the strategy of proofs are ours.


\begin{thebibliography}{99}

	\bibitem{ArcDomPet2016a}
	N.~Arcozzi, K.~Domelevo, and S.~Petermichl.
	\newblock Second {O}rder {R}iesz {T}ransforms on {M}ultiply--{C}onnected {L}ie
		{G}roups and {P}rocesses with {J}umps.
	\newblock {\em Potential Anal.}, 45(4):777--794, 2016.

	\bibitem{BW1995}
	R.~Ba\~nuelos and G.~Wang,
	Sharp inequalities for martingales with applications to the Beurling--Ahlfors and Riesz transforms,
	\emph{Duke Math.\ J.}~\textbf{80} (1995), 575--600.

	\bibitem{BJ2008}
	R.~Ba\~nuelos and P.~Janakiraman,
	$L^p$-bounds for the Beurling--Ahlfors transform,
	\emph{Trans.\ Amer.\ Math.\ Soc.}~\textbf{360} (2008), 3603--3612.

	\bibitem{BK2019}
	R.~Ba\~nuelos and M.~Kwa\'snicki,
	On the~$\ell^p$-norm of the discrete Hilbert transform,
	\emph{Duke Math.\ J.}~\textbf{168} (2019), 471--504.

	\bibitem{Laeng2007}
	E.~Laeng,
	Remarks on the Hilbert transform and some families of multiplier operators related to it,
	\emph{Collect.\ Math.}~\textbf{58} (2007), 25--44.

	\bibitem{BMH2003}
	R.~Ba\~nuelos and P.~J.~M\'endez-Hern\'andez,
	Space-time Brownian motion and the Beurling--Ahlfors transform,
	\emph{Indiana Univ.\ Math.\ J.}~\textbf{52} (2003), 981--990.

	\bibitem{BJV2011}
	A.~Borichev, P.~Janakiraman and A.~Volberg,
	Subordination by orthogonal martingales in~$L^p$ and zeros of Laguerre polynomials,
	preprint, arXiv:1012.1135; see also related estimates in
	\emph{Adv.\ Math.}~\textbf{226} (2011).

	\bibitem{Burkholder1979}
	D.~L.~Burkholder,
	A sharp inequality for martingale transforms,
	\emph{Ann.\ Probab.}~\textbf{7} (1979), 858--863.

	\bibitem{Burkholder1984}
	D.~L.~Burkholder,
	Boundary value problems and sharp inequalities for martingale transforms,
	\emph{Ann.\ Probab.}~\textbf{12} (1984), 647--702.

	\bibitem{Burkholder1991}
	D.~L.~Burkholder,
	Explorations in martingale theory and its applications,
	in \emph{\'Ecole d'\'Et\'e de Probabilit\'es de Saint-Flour XIX---1989},
	Lecture Notes in Math.~\textbf{1464}, Springer, 1991, pp.~1--66.

	\bibitem{DIPV2026}
	K.~Domelevo, P.~Ivanisvili, S.~Petermichl and A.~Volberg,
	Dimension-free bounds for Riesz transforms on the Hamming cube via a Bellman function,
	preprint, arXiv:2606.20289, 2026.

	\bibitem{DP2014}
	K.~Domelevo and S.~Petermichl,
	Sharp~$L^p$ estimates for discrete second-order Riesz transforms,
	\emph{Adv.\ Math.}~\textbf{262} (2014), 932--952.

	\bibitem{Essen1984}
	M.~Ess\'en,
	A superharmonic proof of the M.~Riesz conjugate function theorem,
	\emph{Ark.\ Mat.}~\textbf{22} (1984), 241--249.

	\bibitem{GMSS2009}
	S.~Geiss, S.~Montgomery-Smith and E.~Saksman,
	On singular integral and martingale transforms,
	\emph{Trans.\ Amer.\ Math.\ Soc.}~\textbf{362} (2010), 553--575.

	\bibitem{GV1976}
	R.~F.~Gundy and N.~Th.~Varopoulos,
	Les translations des martingales et les int\'egrales stochastiques,
	\emph{C.\ R.\ Acad.\ Sci.\ Paris S\'er.\ A-B}~\textbf{283} (1976), A509--A511.

	\bibitem{JMP2018}
	M.~Junge, T.~Mei and J.~Parcet,
	Noncommutative Riesz transforms---dimension free bounds and Fourier multipliers,
	\emph{J.\ Eur.\ Math.\ Soc.}~\textbf{20} (2018), 529--584.

	\bibitem{LP2004}
	F.~Lust-Piquard,
	Dimension free estimates for discrete Riesz transforms on products of abelian groups,
	\emph{Adv.\ Math.}~\textbf{185} (2004), 289--327.

	\bibitem{NazTre1996a}
	F.~Nazarov and S.~Treil.
	\newblock The hunt for a {B}ellman function: applications to estimates for
	singular integral operators and to other classical problems of harmonic
	analysis.
	\newblock {\em Algebra i Analiz}, 8(5):32--162, 1996.

	\bibitem{NV2004}
	F.~Nazarov and A.~Volberg,
	Heating of the Ahlfors--Beurling operator, and estimates of its norm,
	\emph{St.\ Petersburg Math.\ J.}~\textbf{15} (2004), 563--573.

	\bibitem{PV2002}
	S.~Petermichl and A.~Volberg,
	Heating of the Ahlfors--Beurling operator: weakly quasiregular maps on the plane are quasiregular,
	\emph{Duke Math.\ J.}~\textbf{112} (2002), 281--305.

	\bibitem{Pichorides1972}
	S.~K.~Pichorides,
	On the best values of the constants in the theorems of M.~Riesz, Zygmund and Kolmogorov,
	\emph{Studia Math.}~\textbf{44} (1972), 165--179.

	\bibitem{Tre2013b}
	S.~Treil.
	\newblock Sharp {$A_2$} estimates of {H}aar shifts via {B}ellman function.
	\newblock In {\em Recent trends in analysis}, volume~16 of {\em Theta Ser.\ Adv.\ Math.},
	pages 187--208. Theta, Bucharest, 2013.

	\bibitem{Wang1995}
	G.~Wang,
	Differential subordination and strong differential subordination for continuous-time martingales and related sharp inequalities,
	\emph{Ann.\ Probab.}~\textbf{23} (1995), 522--551.

\end{thebibliography}
\end{document}